\documentclass[12pt,reqno]{amsart}

\usepackage{amscd}

\let\eps\varepsilon
\def\Type{\par\bgroup\def\endType{\par\egroup}\noindent\tt\raggedright}

\numberwithin{equation}{section}

\theoremstyle{plain}
\newtheorem{thm}{Theorem}
\newtheorem{cor}[thm]{Corollary}
\newtheorem{lem}[thm]{Lemma}
\newtheorem{prop}[thm]{Proposition}
\newtheorem{conj}[thm]{Conjecture}

\theoremstyle{remark}

\newtheorem{example}{Example}

\def\C{\mathbb C}

\def\Z{\mathbb Z}

\def\Hch{\mathcal H'}
\def\Sym{{\mathop{Sym}}}

\def\Sing{{\operatorname{Sing}}}
\def\sing{{\operatorname{sing}}}

\def\Hom{{\operatorname{Hom}}}
\def\End{{\operatorname{End}}}
\def\Res{{\operatorname{Res}}}
\def\Id{{\operatorname{Id}}}
\def\im{{\operatorname{im \ }}}

\def\wt{{\operatorname{wt}}}

\def\Re{\operatorname{Re}}
\def\Im{\operatorname{Im}}

\def\1{{\mathbf 1}}

\def\a{\alpha}

\def\gl{\mathfrak {gl}}
\def\sl{\mathfrak {sl}}

\def\g{\mathfrak g}
\def\h{\mathfrak h}
\def\n{\mathfrak n}
\def\<{\langle}
\def\>{\rangle}
\def\({\left(}
\def\){\right)}

\begin{document}

\title[How to regularize singular vectors and kill the dynamical Weyl group]
{\strut\\[-3.3\baselineskip]
How to regularize singular vectors\\[4pt] and kill the dynamical Weyl group}
\author[K. Styrkas, V. Tarasov, and A. Varchenko]{K. Styrkas$^{\,\star}$,
V. Tarasov$^{\,\diamond,1}$, and A. Varchenko$^{\,\star,2}$}

\thanks{$^{1\,}$Supported in part
by RFFI grant 02-01-00085a and CRDF grant RM1-2334-MO-02}
\thanks{$^{2\,}$Supported in part by NSF grant DMS-9801582}

\maketitle
\begin{center}
{\it
$^\diamond$St.Petersburg Branch of Steklov Mathematical Institute\\
Fontanka 27, St.Petersburg 191011, Russia

\medskip
$^\star$Department of Mathematics,
University of North Carolina at Chapel Hill\\
Chapel Hill, NC 27599-3250, USA
}
\end{center}

%\date{}

\begin{abstract}
Let $\g$ be a simple Lie algebra, and let $M_\lambda$ be the Verma module over
$\g$ with highest weight $\lambda$. For a finite-dimensional $\g$-module $U$
we introduce a notion of a regularizing operator, acting in $U$, which makes
the meromorphic family of intertwining operators
$\Phi:M_{\lambda+\mu}\to M_\lambda\otimes U$ holomorphic, and conjugates
the dynamical Weyl group operators $A_w(\lambda)\in\End(U)$ to constant
operators. We establish fundamental properties of regularizing operators,
including uniqueness, and prove the existence of a regularizing operator
in the case $\g=\sl_3$.
\end{abstract}

\thispagestyle{empty}

%%%%

\section{Introduction}
Let $\g$ be a simple Lie algebra. Let $M_\lambda$ be the Verma module over
$\g$ with highest weight $\lambda$, and $\1_\lambda\in M_\lambda$ the
highest weight vector. A fundamental object of the representation theory
of $\g$ is a family of intertwining operators
$$
\Phi_\lambda^u\ :\ M_{\lambda+\mu} \ \to\ M_\lambda \otimes U , \qquad \1_{\lambda+\mu}\ \mapsto
\ \1_\lambda \otimes u \ + \ \text{ lower order terms },
$$
where $U$ is a finite-dimensional $\g$-module, and $u \in U[\mu]$.
In this paper we study singularities of $\Phi_\lambda^u$ as a function
of $\lambda$ and $u$.

The intertwining operator $\Phi_\lambda^u$ is completely determined by the
image of the highest weight vector, which we denote by $\Sing(\1_\lambda
\otimes u)$. For any $u \in U[\mu]$, the singular vector $\Sing(\1_\lambda
\otimes u)$ can be written in terms of the Shapovalov form on $M_\lambda$,
and is unique for generic $\lambda$. It is a rational function of $\lambda$
with possible poles if $\lambda$ belongs to at least one of the Kac-Kazhdan
hyperplanes $\<\alpha,\lambda+\rho\>=\frac k2\<\alpha,\alpha\>$. For those
special values of $\lambda$ the operator $\Phi_\lambda^u$ may fail to exist.

One of the objectives of this paper is to regularize the family of operators
$\Phi_\lambda^u$, that is, to construct a family of intertwining operators
$\tilde\Phi_\lambda^u: M_{\lambda+\mu}\ \to\ M_\lambda\otimes U$,
holomorphically depending on $\lambda$ and $u$, such that the correspondence
$u\mapsto\tilde\Phi_\lambda^u$ is still injective.

We also study the relation between the operators $\Phi_\lambda^u$ for different
values of $\lambda$ and $u$. Recall that for a dominant integral weight $\nu$ and any
element $w$ of the Weyl group $W$, the Verma module $M_{\nu}$ contains
a submodule isomorphic to $M_{w\cdot\nu}$. We fix this identification and
regard $M_{w\cdot\nu}$ as a submodule of $M_\nu$.

If $\lambda+\mu$ is a dominant integral weight, then for any $w \in W$ the
restriction of the operator $\Phi_\lambda^u$ to the submodule $M_{w \cdot
(\lambda+\mu)}$ takes values in the submodule ${M_{w\cdot\lambda}\otimes U}$.
The relation between this restriction and the operator $\Phi_{w\cdot\lambda}^u$
is given by
$$
\Phi_\lambda^u |_{M_{w\cdot\lambda}}=\Phi_{w\cdot\lambda}^{A_w(\lambda) u}\,,
$$
where $A_w(\lambda)$ is a rational $\End(U)$-valued function of $\lambda$.
The operators $A_w(\lambda)$, $w\in W$, are called the dynamical Weyl group
operators.

The specialization of the dynamical Weyl group operators at $\lambda = -\rho$
produces constant operators $\bar w = A_w(-\rho)$, and for a simple root
reflection $s_i \in W$ we have $\bar s_i = f_i^{h_i}$. According to the
classical results of Verma and Lusztig, the operators $\bar s_i$ form a
 representation of the Weyl group. The dynamical Weyl group
operators may be regarded as meromorphic deformations of the operators
$\bar w$ in the class of cocycles on the Weyl group, satisfying
$$
A_{w_1w_2}(\lambda) = A_{w_1}(w_2\cdot \lambda) A_{w_2}(\lambda)\,.
$$

We introduce a notion of a regularizing operator $N(\lambda)$, acting in $U$
and polynomially depending on $\lambda$. By definition, we require that
$N(\lambda)$ has the following properties.

First, the modified intertwining operators
$\Phi_\lambda^{N(\lambda)u}:M_{\lambda+\mu}\to M_\lambda\otimes U$ are
holomorphic in $\lambda$ for any $u\in U[\mu]$.

Second, the conjugated dynamical Weyl group operators
$N(w\cdot\lambda)^{-1} A_w(\lambda) N(\lambda)$ are constant operators. This
implies, in particular, that the cocycle $A_w(\lambda)$ is cohomologically
equivalent to a constant function of $\lambda$.

Third, we impose a minimality condition on $N(\lambda)$ by requiring that
$\det N(\lambda)$ is equal to a certain explicitly specified polynomial.
(The second property of $N(\lambda)$ implies the divisibility of
$\det N(\lambda)$ by this polynomial.)

We prove that the regularizing operator $N(\lambda)$, if it exists, is uniquely
determined up to the right multiplication by an operator, polynomially
depending on $\lambda$, invertible for all values of $\lambda$, and symmetric
with respect to the Weyl group action. If $\g=\sl_2$ it is easy to check that
the regularizing operator exists for any irreducible finite-dimensional $\g$-module 
$U$; in each weight subspace of $U$ the operator $N(\lambda)$ acts as a scalar,
determined by the determinant condition. In this paper we prove the existence
property in the case $\g=\sl_3$. We conjecture that a regularizing operator
exists for any simple Lie algebra $\g$ and any finite-dimensional $\g$-module
$U$.

Our construction of the regularizing operator $N(\lambda)$ relies on a certain
realization of finite-dimensional $\sl_3$-modules, based on a special case
of the $(\gl_m,\gl_n)$ duality. We identify weight subspaces $U[\mu]$ of
an $\sl_3$-module $U$ with subspaces of singular vectors in tensor products
$V_{(l_1,0)} \otimes V_{(l_2,0)} \otimes V_{(l_3,0)}$ of finite-dimensional
$\gl_2$-modules. In this realization the action of operators $A_w(\lambda)$ is
identified with the action of suitably renormalized standard rational
$R$-matrices, the spectral parameters in the $R$-matrices being determined by
$\lambda$.

We use the functional realization $\mathcal H[\mathbf l]$ of the tensor product
$V_{(l_1,0)} \otimes V_{(l_2,0)} \otimes V_{(l_3,0)}$, provided by the
representation theory of the Yangian $Y(\gl_2)$, to define operators $\mathcal
N[\mathbf l; \mathbf z]:\mathcal H[\mathbf l] \to V_{(l_1,0)} \otimes
V_{(l_2,0)} \otimes V_{(l_3,0)}$, which conjugate the renormalized $R$-matrices
to the identity. We then use results on the reducibility of the tensor product
$V_{(l_1,0)}(z_1) \otimes V_{(l_2,0)}(z_2) \otimes V_{(l_3,0)}(z_3)$ of the
evaluation $Y(\gl_2)$-modules, to prove that for any $u \in \operatorname{Im}\,
\mathcal N[\mathbf l; \mathbf z]$ the operator $\Phi_\lambda^u$ is holomorphic
in $\lambda$.

The final step of the construction is an alternative `constant' identification
of the functional space $\mathcal H[\mathbf l]^\sing$ with the subspace of
singular vectors $\(V_{(l_1,0)}\otimes V_{(l_2,0)}\otimes V_{(l_3,0)}\)^\sing$,
which is provided by a certain combinatorial lemma. The composition of the
`constant' identification with $\mathcal N[\mathbf l; \mathbf z]$ gives the
desired regularizing operator $N(\lambda)$ in our particular realization.

As an application of our results, we use the regularized family of intertwining
operators to introduce an $\End(U[0])$-valued function $\Psi(\lambda,x)$ of
$\lambda\in\h^*, x\in\h$, defined as a certain matrix trace of
$\Phi_\lambda^u$. This Baker-Akhiezer type function is holomorphic in
$\lambda$, and satisfies certain algebraic identites (resonance conditions),
relating values of $\Psi(\lambda,x)$ for different values of $\lambda$.
The fact that $N(\lambda)$ conjugates
the dynamical Weyl group operators $A_w(\lambda)$ to the constant operators
$\bar w$ implies that these identities have a simple form
$$
\Psi(\lambda,x)u=\Psi(s_\alpha \cdot \lambda,x)u\,,
$$
for suitable $\lambda,u$. Special cases of such resonance conditions were used
in \cite{ES} to establish the algebraic integrability of the generalized
quantum Calogero-Sutherland systems.

This paper is organized as follows.
In Section 2 we review facts on $\g$-modules and set the notation.
In Section 3 we recall the construction of intertwining operators
$\Phi:M_{\lambda+\mu} \to M_\lambda \otimes U$ and singular vectors in
the tensor product.

In Section 4 we review facts about the Weyl group and the dynamical Weyl group.
We introduce a distinguished action of $W$ in $U$, and show that the dynamical
Weyl group operators can be viewed as its deformation.  We establish Theorems
\ref{thm:dyngroup} and \ref{thm:dynproperty}, which improve results of
\cite{EV,TV}.

In Section 5 we introduce the notion of a regularizing operator $N(\lambda)$,
and use the axiomatic description to establish properties of regularizing
operators. Theorem \ref{thm:uniqueness} gives the uniquiness of a regularizing
operator. We formulate Theorem \ref{thm:main} on the existence of regularizing
operators for $\g=\sl_3$. This is the main technical result of the paper.

Section 6 provides the background necessary to construct regularizing operators
for $\sl_3$. We invoke the $(\gl_m, \gl_n)$-duality and functional spaces
to give a realization of an $\sl_3$-module $U$, and use the realization to
construct the fundamental operator $\mathcal N[\mathbf l; \mathbf z]$,
which is the main ingredient in the construction of $N(\lambda)$.

Section 7 contains the construction of regularizing operators for $\sl_3$
and proofs of all their properties. The short Section 8 is devoted to
resonance conditions for trace functions, renormalized by $N(\lambda)$.

In Appendix A we review the calculus of formal monomials, which is used to
prove the properties of the dynamical Weyl group operators. Appendix B is
devoted to the representation theory of the Yangian $Y(\gl_2)$. We use Yangian
modules to establish properties of the operator
$\mathcal N[\mathbf l;\mathbf z]$. Appendix C contains a proof of a technical
combinatorial Lemma \ref{thm:constantbasis}, which is used in the main
construction.

\section{Notation}

Let $\g$ be a simple Lie algebra over $\C$ with a Cartan subalgebra $\h$ and
a root system $\Delta$ with a polarization $\Delta = \Delta^+ \cup \Delta^-$.
We identify $\h$ with $\h^*$ using the Killing form on $\g$, and denote the
induced invariant bilinear form on $\h^*$ by $\<\cdot,\cdot\>$.

An element $\lambda \in \h^*$ is called dominant, if $\<\alpha,\lambda\> \ge 0$
for any $\alpha \in \Delta^+$. We write $\lambda \ge \mu$ for $\lambda,\mu \in
\h^*$, if $\lambda-\mu$ is dominant. This is a partial order on $\h^*$.

An element $\lambda \in \h^*$ is called integral, if $\frac{2\<\alpha,\lambda\>}{\<\alpha,\alpha\>} \in
\Z$ for any $\alpha \in \Delta^+$.

Denote $W$ the Weyl group of $\g$. We introduce two actions of $W$ in $\h^*$.
The standard action is defined by $$s_\alpha v = v - \frac
{2\<\alpha,v\>}{\<\alpha,\alpha\>} \alpha$$ for any root reflection $s_\alpha
\in W$, and the ``dot'' action is defined by $$s_\alpha \cdot v =
s_\alpha(v+\rho) - \rho,$$ where $\rho = \frac 12 \sum_{\alpha \in \Delta^+}
\alpha$.

Let $\{e_i,h_i,f_i\}$ denote the standard generators of the Lie algebra $\g$, corresponding to a simple root $\alpha_i$. We denote $\sl_2(\alpha_i)$ the Lie subalgebra of $\g$, spanned by $e_i,h_i,f_i$.

We denote $\n^+,\n^-$ the Lie subalgebras, generated respectively by $\{e_i\}$ and $\{f_i\}$. We have a vector space decomposition $\g = \n^+ \oplus \h \oplus \n^-$.

The root lattice $Q$ is defined by $Q = \sum_{i=1}^{\dim \h} \Z \alpha_i$. We also set
$Q^+ = \sum_{i=1}^{\dim \h} \Z_{\ge0} \alpha_i$.

The universal enveloping algebra $\mathcal U(\g)$ is a $Q$-graded associative
algebra, with the grading defined by
$$
\wt (e_i) = \alpha_i, \quad \wt (h_i) =
0, \quad \wt (f_i) = -\alpha_i, \qquad i = 1,\dots,\dim \h.
$$
The Poincare-Birkhoff-Witt theorem gives us a decomposition
$$
\mathcal U(\g) =
\mathcal U(\n^-) \otimes \mathcal U(\h) \otimes \mathcal U(\n^+),
$$
and we denote $\pi_0: \mathcal U(\g) \to \mathcal U(\h)$ the induced projection
along the subspace $\( \n^- \mathcal U(\g) + \mathcal U(\g) \n^+\) \subset
\mathcal U(\g)$.

Introduce an anti-involution $\varpi$ of $\g$ by
$$\varpi(e_i) = f_i, \qquad \varpi(f_i) = e_i, \qquad \varpi(h_i) = h_i, \qquad i = 1,\dots,\dim \h.$$
Define a bilinear $\mathcal U(\h)$-valued form $S$ on $\mathcal U(\g)$, by
$$S(x,y) = \pi_0(\varpi(x) y) \in \mathcal U(\h), \qquad x,y, \in \mathcal U(\g).$$
The form $S$ is contravariant, i.e. satisfies
$$S(x y,z) = S(y,\varpi(x)z), \qquad x,y,z \in \mathcal U(\g).$$
We will identify the commutative algebra $\mathcal U(\h)$ with the algebra $\C[\h^*]$ of polynomial functions on $\h^*$, and for any $\lambda \in \h^*$ we will denote $S_\lambda$ the $\C$-valued contravariant form on $\mathcal U(\g)$, obtained by evaluating $S$ at $\lambda$.

We consider $\g$-modules $V$ with a weight space decomposition
$$V = \bigoplus_{\mu\in\h^*} V[\mu], \qquad V[\mu] = \{ v \in V\quad|\quad h v = \mu(h) v \text{ for all } h \in \h \},$$
with finite-dimensional weight subspaces $V[\mu]$. We say that $v \in V$ is homogeneous of weight $\mu$ and write $\wt(v) = \mu$, if $v \in V[\mu]$. Then for any homogeneous $x \in \mathcal U(\g)$, we have
$\wt(x v) = \wt (x) + \wt (v)$.

In this paper $U$ always denotes a finite-dimensional $\g$-module.

For any $\lambda \in \h^*$, we denote $M_\lambda$ the Verma module for $\g$,
generated by the highest weight vector $\1_\lambda$, with relations
$$\n^+\1_\lambda = 0, \qquad h \1_\lambda = \lambda(h) \1_\lambda,
\ \text{ for } h \in
\h.$$ We identify any Verma module $M_\lambda$ with $\mathcal U(\n^-)$ by
sending $x \in \mathcal U(\n^-)$ to $x \1_\lambda \in M_\lambda$.

The Verma module $M_\lambda$ carries a bilinear $\g$-contravariant form
$(\cdot|\cdot)$, defined by $$(x \1_\lambda| y \1_\lambda) =
S_\lambda(x,y),\qquad x,y \in \mathcal U(\g).$$ This form is called the
Shapovalov form. Different weight subspaces of $M_\lambda$ are orthogonal with
respect to this form, and its restriction to a weight subspace $M_\lambda[\mu]$
is nondegenerate unless $\lambda$ belongs to a hyperplane of the form
$$\chi_{\alpha,k}(\lambda) \stackrel{def}=
\frac{2\<\a,\lambda+\rho\>}{\<\alpha,\alpha\>} - k = 0$$ for some $\alpha \in
\Delta^+, k \in \Z_{\ge0}$.

The Verma module $M_\lambda$ is reducible if and only if $\lambda$ belongs
to at least one of those hyperplanes. The kernel of the Shapovalov form is
the maximal proper submodule of $M_\lambda$. Let $V_\lambda$ denotes
the irreducible quotient of the Verma module $M_\lambda$; it inherits
a contravariant form and a weight decomposition from $M_\lambda$; moreover,
the contraviriant form on $V_\lambda$ is nondegenerate. It is known that
the module $V_\lambda$ is finite-dimensional if and only if $\lambda$ is
a dominant integral weight.

%%%%

\section{Singular vectors and intertwining operators}

A vector $v$ in a $\g$-module $V$ is called singular, if $\n^+ v=0$.
The subspace of all singular vectors in $V$ is denoted $V^\sing$.
For any $\lambda\in\h^*$, we have an isomorphism
$$
\Hom_\g(M_\lambda,V)\cong V^\sing[\lambda],
$$
constructed by associating with any intertwining operator
$\Phi\in\Hom_\g(M_\lambda,V)$ the image of the vector $\1_\lambda$.

It is known \cite{BGG} that a Verma module $M_\lambda$ contains a unique up to proportionality singular vector of weight $\nu$ if and only if there exists a finite sequence of weights
$$\lambda^{(0)} = \lambda, \lambda^{(1)}, \dots, \lambda^{(k)} = \nu,$$
such that
$$\lambda^{(i)} = s_{\beta_i} \cdot \lambda^{(i-1)} = \lambda^{(i-1)} - n_i \beta_i, \qquad i =1,\dots,k,$$
for some $\beta_i \in \Delta^+, n_i \in \Z_{\ge0}$.
We write $\nu \prec \lambda$ if this condition is satisfied.

\begin{thm} \label{thm:ffm} \cite{FFM}
Let $\lambda \in \h^*$, and let $w \in W$ be such that $w\cdot \lambda \prec
\lambda$. Let $w = s_{i_1} \dots s_{i_l}$ be the reduced decomposition of $w$
in terms of simple root reflections. Then the formal monomial
$$
F_w(\lambda) =
f_{i_1}^{\gamma_1(\lambda)} \dots f_{i_l}^{\gamma_l(\lambda)}, \qquad
\gamma_k(\lambda) = \<\alpha_{i_k},(s_{k+1}\dots s_l) \cdot \lambda\>,
$$
makes sense, and $F_w(\lambda) \1_\lambda$ is a singular vector of weight
$w\cdot\lambda$ in the Verma module $M_\lambda$.
\end{thm}

If $\lambda$ is a dominant integral weight, then $w\cdot\lambda\prec\lambda$
for any $w\in W$, and $\gamma_k(\lambda) \in \Z_{\ge0}$ for all $k=1,\dots,l$.
Therefore, $F_w(\lambda)$ is a well-defined element in $\mathcal U(\n^-)$.
The meaning of $F_w(\lambda)$ for other $\lambda \in \h^*$, when we may have
$\gamma_k(\lambda) \notin \Z_{\ge_0}$, is explained in Appendix A.

Introduce notation $v_{w \cdot \lambda} = F_w(\lambda) \1_\lambda$ for the singular vector, constructed in Theorem \ref{thm:ffm}. We identify the submodule of $M_\lambda$, generated by $v_{w \cdot \lambda}$, with $M_{w\cdot\lambda}$ by mapping $\1_{w\cdot\lambda}$ to $v_{w\cdot\lambda}$.

We now study singular vectors in tensor products $M_\lambda\otimes U$,
$\lambda\in\h^*$. All tensor products in this paper are considered over
the field $\C$.

Let $\{g_i\}$ be a homogeneous basis in $\mathcal U(\n^-)$, and let
$\(S_\lambda^{-1}\)_{ij}$ be the inverse to the matrix of the form $S_\lambda$,
restricted to $\mathcal U(\n^-)$, with respect to this basis. The matrix
elements $\(S_\lambda^{-1}\)_{ij}$ are rational functions of $\lambda\in\h^*$,
with possible simple poles at the hyperplanes $\chi_{\alpha,k}(\lambda) = 0$,
for some $\a\in\Delta^+$ and $k\in\Z_{\ge0}$.

Introduce the Cartan anti-automorphism $\omega$ of the Lie algebra $\g$ by
$$
\omega(e_i) = -f_i, \qquad \omega(f_i) = -e_i, \qquad \omega(h_i) = h_i,
\qquad i = 1,\dots,\dim \h\,.
$$
Consider a rational $\mathcal U(\n^-) \otimes \End_\C(U)$-valued function of
$\lambda \in \h^*$,
\begin{equation}
\Xi(\lambda) = \sum_{i,j}
\(S_\lambda^{-1}\)_{ij} g_i \otimes \omega(g_j)\,,
\label{eq:Xi}
\end{equation}
where we regard $\omega(g_j)$ as operators acting in $U$. Note that the
summation over $i,j$ is finite, because $\(S_\lambda^{-1}\)_{ij} = 0$ for
sufficiently large $i$ and $\omega(g_j)$ acts as zero for sufficiently large
$j$. One can check that $\Xi(\lambda)$ does not depend on the choice of the
basis $\{g_j\}$.

For any $\lambda\in\h^*$ and $u\in U$ such that $\Xi(\lambda) (1 \otimes u) \in
\mathcal U(\n^-) \otimes U$ is well-defined, we denote $$\Sing(\1_\lambda
\otimes u) = \Xi(\lambda) (\1_\lambda \otimes u) \in M_\lambda \otimes U.$$ For
generic $\lambda$, the vector $\Sing(\1_\lambda \otimes u)$ is the unique
singular vector of the form
$$
\Sing(\1_\lambda \otimes u) = \1_\lambda \otimes
u + \text{ lower order terms },
$$
where the ``lower order terms'' is a linear combination of vectors $v'\otimes
u'$ with $\wt(v') < \lambda$, see \cite{ES}. The vector $\Sing(\1_\lambda
\otimes u)$ is a singular vector for any $\lambda$, when it is well-defined.

For any $u\in U[\mu]$ such that $\Sing(\1_\lambda\otimes u)$ is well defined
we introduce an intertwining operator $\Phi_\lambda^u \in \Hom_\g
(M_{\lambda+\mu},M_{\lambda} \otimes U)$, uniquely determined by the condition
$$
\Phi_\lambda^u \1_{\lambda+\mu} = \Sing(\1_\lambda \otimes u).
$$

\section{Weyl group and dynamical Weyl group}

The Weyl group $W$, associated with $\g$, is generated by simple reflections
$\{s_i\},\ i=1,\dots,\dim \h$, subject to the relations $s_i^2 = id$,
and the braid relations
\begin{equation}
\underbrace{s_i s_j s_i \dots}_{n_{ij} \text{ factors }} = \underbrace{s_j s_i s_j \dots}_{n_{ij} \text{ factors }}, \label{eq:braid}
\end{equation}
for $i \ne j$, where $n_{ij} = 2,3,4,6$ if
$\frac{4\<\alpha_i,\alpha_j\>^2}{\<\alpha_i,\alpha_i\>\<\alpha_j,\alpha_j\>} =
0,1,2,3$ respectively.

The length $l(w)$ of any $w \in W$ is defined to be the smallest positive
integer $l$ such that $w$ can be written as $w = s_{i_1} \dots s_{i_l}$ for
some $i_1,\dots,i_l$. For the identity element $id \in W$ we set $l(id) = 0$.

Let $U$ be a finite-dimensional $\g$-module. For any $i = 1,\dots,\dim \h$, introduce an operator
$\bar s_i \in \End_\C(U)$ by setting
\begin{equation}
\bar s_i f_i^k v = f_i^{n-k} v, \qquad k = 0,1, \dots, n, \label{eq:Weylaction}
\end{equation}
for any $v \in U$ such that $e_i v = 0$ and $h_i v = n v$ for some $n \in \Z_{\ge0}$.
One can think of operators $\bar s_i$ as operators $f_i^{h_i}$.

The following lemma is well-known; see \cite{EV} and references therein.

\begin{lem}
\label{thm:wbarep}
There exists a linear representation $w \mapsto \bar w$ of the Weyl group $W$
in $U$, such that the generators $s_i$ are mapped to the operators $\bar s_i$.
\end{lem}

Obviously, for any $w \in W$ we have
$$
\bar w = \bigoplus \bar w_\mu, \qquad \bar w_\mu \in \Hom_\C(U[\mu],U[w\mu]).
$$
Note that if $w \in W$ and $\mu \in \h^*$ are such that $w\mu = \mu$, then
$\bar w_\mu = \Id_{U[\mu]}$. One can show that the equivalence class of
the representation $w\mapsto \bar w$ of $W$ is completely determined by these
conditions.

{\bf Remark.} The operators $\bar w$ do not preserve the Shapovalov form on $U$;
one might try to modify the formula \eqref{eq:Weylaction}, and introduce
operators $\tilde s_i$ by
\begin{equation}
\tilde s_i f_i^k v = \frac {k!}{(n-k)!} f_i^{n-k} v, \qquad k = 0,1, \dots, n,
\end{equation}
for any $v \in U$ such that $e_i v = 0$ and $h_i v = n v$ for some $n \in \Z_{\ge0}$.
The new operators $\tilde s_i$ clearly preserve the Shapovalov form.

\begin{conj}\label{thm:Weylconj}
There exists a representation of the Weyl group $W$ in $U$, such that the
generators $s_i$ are mapped to the operators $\tilde s_i$.
\end{conj}

This conjecture can be verified in certain cases, including adjoint
representations for arbitrary $\g$; however, no general proof is available.

We now study the dynamical Weyl group, introduced in \cite{TV}, \cite{EV}.
The following Theorem \ref{thm:dyngroup} and Theorem \ref{thm:dynproperty}
 are refinements of previous results.

Let $\C[\h^*], \C(\h^*)$ denote respectively the algebras of polynomial and
rational functions on $\h^*$. If $V$ is a vector space and $F \in V \otimes
\C[\h^*]$, we denote $F(\lambda)\in V$ the value of the $V$-valued polynomial
$F$ at point $\lambda \in \h^*$.

\begin{thm}\label{thm:dyngroup}
For any $w \in W$ there exists an operator $A_w\in \End_\C(U) \otimes \C(\h^*)$, satisfying
the following property.

Suppose $\lambda,\mu \in \h^*$ are such that $w \cdot (\lambda+\mu)\prec \lambda+\mu$, and $u \in U[\mu]$ is such that $\Sing(\1_\lambda \otimes u)$ is well-defined. Then
$A_w(\lambda)u$ is well-defined, and
\begin{enumerate}
\item
If $w \cdot \lambda \prec \lambda$, then
\begin{equation}
\Phi_\lambda^u v_{w\cdot(\lambda+\mu)} = v_{w\cdot\lambda} \otimes A_w(\lambda)u + \text{ lower order terms }. \label{eq:dyngroup}
\end{equation}
\item
If the condition $w\cdot\lambda \prec \lambda$ fails, then $A_w(\lambda)u = 0$.
\end{enumerate}
\end{thm}

This theorem is proved in Appendix A.

The collection of operators $A_w \in \End_\C(U)\otimes \C(\h^*), w \in W$, is called the dynamical Weyl group.
Equation \eqref{eq:dyngroup} implies that
$$A_w = \bigoplus_\mu A_w^\mu, \qquad A_w^\mu \in \Hom_\C(U[\mu],U[w\mu]) \otimes \C(\h^*).$$

{\bf Example.}
Let $U$ be a finite-dimensional $\sl_2$-module, generated by a highest weight
vector $v$ of weight $n \in \Z_{\ge0}$. Then the dynamical Weyl group consists
of two operators: $A_{id} \equiv \Id$ and $A_{s_1}$, determined by
\begin{equation}
A_{s_1}(\lambda) f^k v = (-1)^n
\frac{(-t-2)(-t-3)\dots(-t-n+k-1)}{t(t-1)\dots(t-k+1)}
\frac {k!}{(n-k)!} f^{n-k} v, \qquad k = 0,1, \dots, n\,,
\label{eq:dynsl2}
\end{equation}
where $t=(\lambda,\alpha_1)$.

For arbitrary $\g$, the operators $A_{s_i}(\lambda)$ act in any irreducible
$\sl_2(\alpha_i)$ submodule of $U$ by formula \eqref{eq:dynsl2}, with
$t=\frac{2\<\lambda,\alpha_i\>}{\<\alpha_i,\alpha_i\>}$.

\begin{thm}\label{thm:dynproperty}
The dynamical Weyl group operators $A_w, \ w \in W$, acting in a finite-dimensional $\g$-module $U$, have the following properties.

\begin{enumerate}
\item
The dynamical Weyl group operators satisfy the cocycle condition,
\begin{equation}
A_{w_1w_2}(\lambda) = A_{w_1}(w_2\cdot\lambda) A_{w_2}(\lambda), \qquad w_1,w_2 \in W.
\label{eq:cocycle}
\end{equation}
\item
The operator $A_w$, regarded as an $\End_\C(U)$-valued rational function on $\h^*$, has at most simple poles, which may occur only at the hyperplanes
$$\chi_{\alpha,k}(\lambda) = 0, \qquad \alpha \in \Delta^+ \cap w(\Delta^-), \ k \in \Z_{\ge0}.$$
\item
If we identify $U[\mu]$ and $U[w \mu]$ in any way not depending on $\lambda$,
then for generic $\lambda \in\h^*$ we have
$$
\det A_w^\mu(\lambda) = const \ \prod_{\a\in\Delta^+ \cap w(\Delta^-)} \frac
{\prod_{k=0}^\infty \chi_{\alpha,k}(w\cdot\lambda)^{\dim U[w\mu+k\alpha]}}
{\prod_{k=0}^\infty \chi_{\alpha,k}(\lambda)^{\dim U[\mu+k\alpha]}}.
$$
where the constant depends on the choice of identification.
\end{enumerate}
\end{thm}

{\bf Remark.}
The dynamical Weyl group operators $A_w(\lambda)$ were introduced in \cite{TV}, \cite{EV}, using a different normalization of the singular vectors:
$$v_{w\cdot\lambda} = \frac{ f_{i_1}^{n_1}} {n_1!} \dots \frac {f_{i_l}^{n_l}}{n_l!} \1_\lambda, \qquad
n_k = \<\alpha_{i_k},(s_{k+1}\dots s_l) \cdot \lambda\>.$$
Equation \eqref{eq:dyngroup} was observed in \cite{TV}, \cite{EV} only for $\lambda \gg \mu$. Part 2 of Theorem \ref{thm:dyngroup} is new.

In \cite{TV}, \cite{EV}, the cocycle condition was satisfied only for $w_1,w_2
\in W$ such that $l(w_1 w_2) = l(w_1) + l(w_2)$; equivalently, operators
$A_w(\lambda)$ formed a cocycle on the braid group. Our operators satisfy
the additional conditions $A_{s_\alpha} (s_\alpha \cdot \lambda)
A_{s_\alpha}(\lambda) = \Id$, and thus form a cocycle on the Weyl group.

\begin{proof}[Proof of Theorem \ref{thm:dynproperty}]
It follows directly from the construction (see also \cite{TV,EV}) that the cocycle condition \eqref{eq:cocycle} is satisfied when
$w_1,w_2 \in W$ are such that $l(w_1 w_2) = l(w_1) + l(w_2)$. Therefore, it suffices to show that
$$A_{s_i} (s_i \cdot \lambda) A_{s_i}(\lambda) = \Id, \qquad i=1,\dots,\dim \h.$$
This follows from the explicit formula \eqref{eq:dynsl2} for the action of $A_{s_i}(\lambda)$ in any irreducible $\sl_2(\alpha_i)$-submodule of $U$.

Now, let $w\in W$ have a reduced decomposition $w = s_{i_1} \dots s_{i_l}$. Then the cocycle condition implies
$$A_w(\lambda) = A_{s_{i_1}}(s_{i_2}\dots s_{i_l} \cdot \lambda) \dots A_{s_{i_l}}(\lambda),$$
and singularities of $A_w(\lambda)$ are determined by the singularities of $A_{s_{i_k}}(s_{i_{k+1}}\dots s_{i_l} \cdot \lambda)$ for $k=1,\dots,\dim \h$. From the explicit formula \eqref{eq:dynsl2} it follows that
$A_{s_{i_k}}(s_{i_{k+1}}\dots s_{i_l} \cdot \lambda)$ may only have simple poles at the hyperplanes
$\chi_{\alpha_{i_k},n}( s_{i_{k+1}}\dots s_{i_l} \cdot \lambda) = 0$, or equivalently
\begin{equation}
\chi_{\beta_k,n}( \lambda) = 0, \qquad \beta_k = s_{i_l} \dots s_{i_{k+1}}\alpha_{i_k}.
\label{eq:dynpoles}
\end{equation}
It is known that if $w = s_{i_1} \dots s_{i_l}$ is a reduced expression, then
$\{\beta_k\}$ form (without repetitions) the set ${\Delta^+\cap w(\Delta^-)}$.
Therefore, the hyperplanes \eqref{eq:dynpoles} are all distinct, and the second
assertion follows.

Finally, from the cocycle property we obtain
$$\det A_w^\mu(\lambda) = \det A_{s_{i_1}}(s_{i_2}\dots s_{i_l} \cdot \lambda) \dots \det A_{s_{i_l}}(\lambda) $$
and the proof of the general formula reduces to verification of the formula
$$\det A_{s_i}^\mu(\lambda) = const \ \ \frac
{\prod_{k=0}^\infty \chi_{\alpha_i,k}(s_i\cdot\lambda)^{\dim U[s_i\mu+k\alpha_i]}}
{\prod_{k=0}^\infty \chi_{\alpha_i,k}(\lambda)^{\dim U[\mu+k\alpha_i]}}.$$
which follows from \eqref{eq:dynsl2}.
\end{proof}

The cocycle condition \eqref{eq:cocycle} immediately implies that the operators $A_w(-\rho), w \in W$, form a representation of the Weyl group. One can easily check that in fact $A_w(-\rho) = \bar w$, where operators $\bar w$ are defined by \eqref{eq:Weylaction}.

\section{Regularizing operators and their properties}

In this section we define the main object of our study --- regularizing
operators for finite-dimensional representations of a semisimple Lie algebra
$\g$, and establish some of their properties.

Let $\g$ be a simple Lie algebra, and let $U$ be a finite-dimensional $\g$-module.
An operator $N \in \End_\C(U) \otimes \C[\h^*]$ is called a regularizing operator, if it satisfies the following conditions.

\begin{enumerate}
\item
The operator $N$ preserves weight subspaces, i.e. it can be decomposed as
$$N = \bigoplus_\mu N_\mu, \qquad N_\mu \in \End_\C(U[\mu])\otimes \C[\h^*].$$

\item
For any $\mu \in \h^*$, there exists a nonzero constant $c_\mu$ such that
\begin{equation}
\det N_\mu(\lambda) = c_\mu \prod_{\a\in\Delta^+} \prod_{k=0}^\infty
\chi_{\alpha,k}(\lambda)^{\dim U[\mu+k\alpha]}. \label{eq:detN}
\end{equation}

\item
There exist operators $a_w \in \End_\C(U)$, such that for any $w \in W$ and generic $\lambda\in \h^*$ we have
\begin{equation}
a_w = N(w\cdot\lambda)^{-1} A_w(\lambda) N(\lambda). \label{eq:dynN}
\end{equation}

\item \label{eq:regularityN}
The operator
$\Xi_N(\lambda): U \to \mathcal U(\n^-) \otimes U$, defined by
$$
\Xi_N(\lambda)u = \Xi(\lambda) (1 \otimes N(\lambda)u) , \qquad u \in U,
$$
depends polynomially on $\lambda$. Here $\Xi(\lambda)$ is given
by \eqref{eq:Xi}.
\end{enumerate}

{\bf Remark.}
The last condition is equivalent to the requirement that for any $u \in U$, the singular vector
$\Sing(\1_\lambda \otimes N(\lambda)u) \in M_\lambda \otimes U$ is polynomial in $\lambda$.

The following Lemma is convenient for establishing the formula \eqref{eq:detN}.

\begin{lem}\label{thm:detfrompoles}
Suppose $N \in \End_\C(U) \otimes \C[\h^*]$ satisfies conditions (1) and (3) above, and is such that
$N^{-1}$ is regular outside the hyperplanes $\chi_{\alpha,k}(\lambda)=0$ for all $\alpha \in \Delta^+$ and $k \in \Z_{\ge0}$. Then $N$ also satisfies condition (2).
\end{lem}

\begin{proof}
The formula \eqref{eq:dynN} implies that for any $w \in W$ and generic $\lambda \in \h^*$, we have
\begin{equation}
\det A_w^\mu(\lambda) \det N_\mu(\lambda) = const \ \det N_{w\mu}(w\cdot\lambda).
\label{eq:divN}
\end{equation}
The right-hand side of this equation polynomially depends on $\lambda$, and
therefore in the left-hand side $\det N_\mu(\lambda)$ must be divisible by the
denominator of $\det A_w^\mu(\lambda)$. In the special case when $w=w_0$ is the
longest element of the Weyl group, this denominator is precisely
$\prod_{\a\in\Delta^+} \prod_{k=0}^\infty \chi_{\alpha,k}(\lambda)^{\dim
U[\mu+k\alpha]}$. Therefore, we conclude that
$$
c_\mu(\lambda) = \frac {\det N_\mu(\lambda)} {\prod_{\a\in\Delta^+}
\prod_{k=0}^\infty \chi_{\alpha,k}(\lambda)^{\dim U[\mu+k\alpha]}}
$$
is polynomial in $\lambda$, and \eqref{eq:divN} implies that
$c_\mu(\lambda) = const \ c_{w_0\mu}(w_0\cdot\lambda)$.

By the assumption of the Lemma, $c_\mu(\lambda)$ may only vanish at hyperplanes
$\chi_{\alpha,k}(\lambda) = 0$, and similarly $c_{w_\mu}(w_0 \cdot \lambda)$
may only vanish at hyperplanes $\chi_{\alpha,k}(w_0\cdot \lambda) = 0$. Since
these two families of hyperplanes are disjoint, we conclude that
$c_\mu(\lambda)$ is a polynomial in $\lambda$ which never vanishes, i.e. a
constant polynomial. The condition (2) follows.
\end{proof}

The argument above also shows that the condition (2) can be thought of
as a minimality condition for $N$.

Here are some other immediate corollaries of the definition.

\begin{prop}
Let $U$ be a finite-dimensional $\g$-module, and let $N$ be a regularizing operator.

\begin{enumerate}
\item
The operators $a_w \in \End_\C(U)$ form a representation of the Weyl group.

\item
The inverse operator $N^{-1}$, regarded as an $\End_\C(U)$-valued rational function on $\h^*$, may have only simple poles. The poles may only occur at hyperplanes
$\chi_{\alpha,k}(\lambda) = 0$ for $\alpha \in \Delta^+$ and $k \in \Z_{\ge0}$
such that $U[\mu+k\alpha] \ne 0$.

\item
Let $C \in \End_\C(U)$ be a weight-preserving operator, i.e.
$$C = \bigoplus_\mu C_\mu, \qquad C_\mu \in \End(U[\mu]). $$
Then the operator $\bar N = N C$ is also a regularizing operator.

\end{enumerate}
\end{prop}

\begin{proof}
For any $w_1,w_2 \in W$, we compute
\begin{gather*}
a_{w_1w_2} = N(w_1w_2\cdot \lambda)^{-1} A_{w_1w_2}(\lambda) N(\lambda) =
N(w_1w_2\cdot \lambda)^{-1} A_{w_1}(w_2\cdot\lambda) A_{w_2}(\lambda) N(\lambda) = \\
\( N(w_1w_2\cdot \lambda)^{-1} A_{w_1}(w_2\cdot\lambda) N(w_2\cdot\lambda) \)
\( N(w_2\cdot \lambda)^{-1} A_{w_2}(\lambda) N(\lambda) \) = a_{w_1} a_{w_2},
\end{gather*}
which shows that the correspondence $w \mapsto a_w$ gives a representation of the Weyl group.

{}From the formula for $\det N(\lambda)$, we see that $N^{-1}$ may only have
poles at hyperplanes $\chi_{\alpha,k}(\lambda) = 0$, for some
$\alpha\in\Delta^+$ and $k \in \Z_{\ge0}$ such that $U[\mu+k\alpha] \ne 0$.
Fix $\alpha\in\Delta_+$ and $k\in\Z_{\ge0}$. We have
$$
N(\lambda)^{-1} =
a_{s_\alpha}^{-1} N(s_\alpha\cdot\lambda)^{-1} A_{s_\alpha}(\lambda)\,.
$$
since $N(s_\alpha\cdot\lambda)^{-1}$ does not have a pole at the hyperplane
$\chi_{\alpha,k}(\lambda) = 0$, the singularity of $N(\lambda)^{-1}$ may come
only from the operator $A_{s_\alpha}(\lambda)$, which have at most simple pole
there.

The operators $\bar N(\lambda) = N(\lambda)C$ are obviously polynomial in $\lambda$, and the conditions on $C_\mu$ ensure weight preserving and determinant properties of $\bar N(\lambda)$. Finally,
$$\bar a_w = \bar N(w\cdot \lambda)^{-1} A_w(\lambda) \bar N(\lambda) = C^{-1} N(w\cdot \lambda)^{-1} A_w(\lambda) N(\lambda) C = C^{-1} a_w C$$
do not depend on $\lambda$, and for any $u \in U$ the vector
$$\Xi_{\bar N}(\lambda) u = \Xi(\lambda) (1 \otimes \bar N(\lambda)u) = \Xi(1 \otimes N(\lambda) C u) = \Xi_N(\lambda)(Cu)$$
is polynomial in $\lambda$.
\end{proof}

\begin{thm}
Let $U$ be a finite-dimensional $\g$-module, and let $N \in \End_\C(U) \otimes \C[\h^*]$ be a regularizing operator. Fix any $\lambda_0 \in \h^*$. Then

\begin{enumerate}
\item
The linear map
$$
\Xi_N(\lambda_0): U \to \(M_{\lambda_0} \otimes U\)^\sing,
\qquad u \mapsto \Sing(\1_{\lambda_0}\otimes N(\lambda_0)u)
$$
is injective.
\item
We have the following description of the image of $N_\mu(\lambda_0)$:
$$
\operatorname{Im} N(\lambda_0) = \{\,u \in U\ \,|\ \,
\Sing(\1_{\lambda_0} \otimes u) \text{ is well-defined}\,\}\,.
$$
\end{enumerate}
\end{thm}

\begin{proof}

Let $\Delta_0$ be the root subsystem of $\Delta$, defined by
$$
\Delta_0 = \{\,\alpha \in \Delta\ |\ \< \alpha, \lambda_0 \> \in \Z\}\,.
$$
We set $\Delta_0^+ = \Delta_0 \cap \Delta^+$, and denote by $W_0$ the subgroup
of $W$, generated by root reflections $s_\alpha, \alpha \in \Delta_0$.

Let $w \in W_0$ be such that $w\cdot\lambda_0$ is antidominant with respect to
$\Delta_0^+$. Then, $\det N(w \cdot \lambda_0) \ne 0$, and
$N(w \cdot\lambda_0)$ is invertible.

According to Theorem \ref{thm:dyngroup},
$$\Phi_{\lambda_0}^{N(\lambda_0) u} v_{w\cdot(\lambda_0+\mu)} = v_{w\cdot\lambda_0} \otimes A_w(\lambda_0)N(\lambda_0) u + \text{ lower order terms }.$$

Suppose for some $u \in U[\mu]$ we have $\Sing(\1_{\lambda_0} \otimes
N(\lambda_0)) u = 0$. Then the intertwining operator
$\Phi_{\lambda_0}^{N(\lambda_0)u}$ is identically zero, and in particular
$$
N(w\cdot \lambda_0) a_w u = A_w(\lambda_0)N(\lambda_0) u = 0.
$$
Since $N(w \cdot\lambda_0)$ and $a_w$ are both invertible, we conclude that
$u=0$, and thus $\Xi_N(\lambda_0)$ is injective.

Now, let $u \in U[\mu]$ be such that $\Sing(\1_{\lambda_0} \otimes u)$ is well-defined. Then the vector $u' = A_w(\lambda_0) u$ is also well-defined, and we have
$$u = A_w(\lambda_0)^{-1} u' = N(\lambda_0) a_w^{-1} N(w\cdot\lambda_0)^{-1} u'
\in \operatorname{Im} N(\lambda_0).$$
Conversely, property \eqref{eq:regularityN} of $N(\lambda_0)$ guarantees that $\Sing(\1_{\lambda_0} \otimes u)$ is well-defined for any $u \in\operatorname{Im} N_\mu(\lambda_0)$.
\end{proof}

Finally, we prove that the regularizing operator is determined uniquely up to right multiplication by a polynomially invertible matrix $C(\lambda)$, satisfying certain $W$-invariance conditions.

\begin{thm}\label{thm:uniqueness}
Let $N,\bar N$ be regularizing operators, and let $a_w, \bar a_w$ be the associated constant representations of the Weyl group. Then the operator
$$C(\lambda) = N(\lambda)^{-1} \bar N(\lambda) \in \End_\C(U)$$
is polynomial in $\lambda$, has a decomposition
$$C(\lambda) = \bigoplus_\mu C_\mu(\lambda), \qquad C_\mu(\lambda) \in \End_\C(U[\mu]),$$
and satisfies
$$C(w \cdot \lambda) \bar a_w = a_w C(\lambda).$$

Moreover, the inverse operator $C(\lambda)^{-1}$ is also polynomial in $\lambda$, and
the Weyl group representations $a_w,\bar a_w$ are equivalent.
\end{thm}

\begin{proof}
Since regularizing operators preserve weight subspaces, we have
$$C_\mu(\lambda) = N_\mu(\lambda)^{-1} \bar N_\mu(\lambda) \in \End_\C(U[\mu]).$$

The only possible singularities of $C(\lambda) = N(\lambda)^{-1} \bar
N(\lambda)$ may be simple poles at one of the hyperplanes $\chi_{\alpha,k}
(\lambda) = 0$. Fix $\alpha \in \Delta^+, k \in\Z_{\ge0}$, and write
$$
N(\lambda)^{-1} = \frac {X(\lambda)}{\chi_{\alpha,k}(\lambda)} + Y(\lambda)
$$
for some $X(\lambda), Y(\lambda) \in \End_\C(U)$, regular at the hyperplane
$\chi_{\alpha,k}(\lambda) = 0$.

It follows from the identity $N(\lambda)^{-1} N(\lambda)=\Id$
that for generic $\lambda_0$ from this hyperplane, we have
$X(\lambda_0) N(\lambda_0) = 0$. Hence,
$$
\ker \operatorname X(\lambda_0) \supset \operatorname{Im} N(\lambda_0)
= \operatorname{Im} \bar N(\lambda_0)\,,
$$
which implies that $C(\lambda)$ is well-defined at $\lambda_0$.
Therefore, $C(\lambda)$ does not have a pole at the hyperplane
$\chi_{\alpha,k} (\lambda) = 0$, which was arbitrary, and
thus $C(\lambda)$ is polynomial in $\lambda$.

In view of \eqref{eq:detN}, we conclude that $\det C_\mu(\lambda) = \( \det
N_\mu(\lambda) \)^{-1} \det \bar N_\mu(\lambda)$ is a nonzero constant.
Hence, each $C_\mu(\lambda)$ is polynomially invertible.

Finally, we compute
$$
C(w \cdot \lambda)^{-1} a_w C(\lambda) = \bar N(w \cdot \lambda)^{-1}
N(w \cdot \lambda) a_w N(\lambda)^{-1} \bar N(\lambda) =
\bar N(w \cdot \lambda)^{-1} A_w(\lambda) \bar N(\lambda) = \bar a_w.
$$
The invariance condition follows. Setting $\lambda = -\rho$, we get
$$
\bar a_w = C(-\rho)^{-1} a_w C(-\rho)
$$
and the representations $a_w, \bar a_w$ of the Weyl group are equivalent.
\end{proof}

{\bf Example:}\enspace
Let $\g=\sl_2$. We identify $\h^*$ with $\C$ by associating $ \h^* \ni
\lambda \leftrightarrow \<\alpha_1,\lambda\> \in \C$. Under this
identification, $\alpha_1\equiv 2$ and $\rho\equiv 1$.

Let $U$ be an irreducible finite-dimensional $\sl_2$-module with highest weight
$\Lambda\in\Z_{\ge 0}$. For $u \in U[\Lambda-2k]$, $k \in \Z_{\ge0}$, we have
$$
\Sing(\1_\lambda \otimes u) = \sum_{j=0}^k
\frac {(-1)^j}{j!\prod_{j=0}^k (\lambda-j)} f^j \1_\lambda \otimes e^j u\,.
$$

The weight subspaces of $U$ are one-dimensional, the operators
$N_{\Lambda-2k}(\lambda)\in\End(U[\Lambda-2k]) \cong \C$ being scalars,
and we take
$$
N_{\Lambda-2k}(\lambda) =
\frac1{k!}\,\prod_{j=0}^{k-1} \chi_{\alpha,j}(\lambda)
= \frac1{k!}\,\prod_{j=0}^{k-1} (\lambda - j)\,.
$$
One can check that the operator $a_s$ acts in $U$ by
$$
a_s : f^k \1_\Lambda \mapsto f^{\Lambda-k}\1_\Lambda,
\qquad k = 0,1,\dots,\Lambda\,.
$$
In particular, if $\Lambda$ is an even integer and $u = f^{\Lambda/2}
\1_\Lambda \in U[0]$, we have $a_s u = u$.

\medskip
We now state our main result on the existence of regularizing operators.

\begin{thm}\label{thm:main}
Let $\g = \sl_3$, and let $U$ be an irreducible finite-dimensional $\g$-module.
Then there exists a regularizing operator $N \in \End_\C(U) \otimes \C(\h^*)$.
Moreover, the associated operators $a_{s_i}$, corresponding to simple root
reflections in $W$, coincide with the operators $\bar s_i$ given by the
formulas \eqref{eq:Weylaction}.

\end{thm}

This theorem will be proved in Section \ref{sec:main}.

{\bf Example.}
Let $U$ be the adjoint representation of $\sl_3$. The zero weight subspace $U[0]$ is the Cartan subalgebra $\h$; the generators $h_1, h_2 \in \sl_3$ form a basis in $U[0]$.

The action of operators $A_w^{[0]}(\lambda)$ on $U[0]$ is completely determined by the matrices
$$
A_{s_1}^{[0]}(\lambda) \begin{pmatrix}
\frac {-\lambda_1-2}{\lambda_1} & \frac {-\lambda_1-1}{\lambda_1} \\
0 & 1
\end{pmatrix},
\qquad
A_{s_2}^{[0]}(\lambda) \begin{pmatrix}
1 & 0 \\
\frac {-\lambda_2-1}{\lambda_2} & \frac {-\lambda_2-2}{\lambda_2}
\end{pmatrix},
$$
where $\lambda_1 = \<\alpha_1,\lambda\>, \lambda_2 = \<\alpha_2,\lambda\>$.
The operator $N_{[0]}(\lambda): U[0] \to U[0]$ is given by the matrix
$$N_{[0]}(\lambda) = \frac 13
\begin{pmatrix}
2\lambda_1^2 + 2\lambda_1 \lambda_2-\lambda_2^2 + 4\lambda_1 - \lambda_2 &
 2\lambda_1^2 + 2 \lambda_1\lambda_2 - \lambda_2^2 + 2\lambda_1 - 2\lambda_2 \\
2\lambda_2^2 + 2\lambda_1 \lambda_2-\lambda_1^2 - 2\lambda_1 + 2\lambda_2 &
 2\lambda_2^2 + 2 \lambda_1\lambda_2 - \lambda_1^2 - \lambda_1 + 4\lambda_2
\end{pmatrix}.
$$

It is a straightforward computation to show that
$$
N_{[0]}(s_1\cdot \lambda)^{-1} A_{s_1}^{[0]}(\lambda) N_{[0]}(\lambda) \equiv
\Id_{U[0]}, \qquad
N_{[0]}(s_2\cdot \lambda)^{-1} A_{s_2}^{[0]}(\lambda) N_{[0]}(\lambda) \equiv
\Id_{U[0]},$$
and that
$$
\det N_{[0]}(\lambda) = \lambda_1 \lambda_2 (\lambda_1+\lambda_2+1) =
\chi_{\alpha_1,1}(\lambda) \chi_{\alpha_2,1}(\lambda)
\chi_{\alpha_1+\alpha_2,1}(\lambda).
$$

%%%%

\section{The $(\gl_m,\gl_n)$ duality,
and functional realizations of $\gl_n$-modules}\label{sec:duality}

In Section \ref{sec:main} we will construct the regularizing operators for $\g = \sl_3$ using a concrete realization of irreducible finite-dimensional $\sl_3$-modules. This section is devoted to a review of some underlying results on representations of the reductive Lie algebras $\gl_n$.

The space $\C^{mn} = \(\C^m\)^{\otimes n} = \(\C^n\)^{\otimes m}$ has natural
commuting actions of $\gl_m$ and $\gl_n$. The algebra $\mathbb P_{m,n} =
S(\C^{mn})$ of polynomial functions on $\C^{mn}$ becomes a module over
$\gl_m\otimes\gl_n$. It has decompositions
$$
\mathbb P_{m,n} =
\bigoplus_{\mu_1,\dots,\mu_m} \mathbb P_{m,n}^{\gl_m}[(\mu_1,\dots,\mu_m)]
= \bigoplus_{\nu_1,\dots,\nu_n} \mathbb P_{m,n}^{\gl_n}[(\nu_1,\dots,\nu_n)]
$$
where the subspaces
\begin{gather*}
\mathbb P_{m,n}^{\gl_m}[(\mu_1,\dots,\mu_m)] \cong
S^{\mu_1}(\C^n) \otimes\dots\otimes S^{\mu_m}(\C^n),
\\
\mathbb P_{m,n}^{\gl_n}[(\nu_1,\dots,\nu_n)] \cong
S^{\nu_1}(\C^m) \otimes\dots\otimes S^{\nu_n}(\C^m),
\end{gather*}
are the weight subspaces with respect to the algebras $\gl_m$ and $\gl_n$
respectively.

In the explicit realization of $\mathbb P_{m,n}$ as the polynomial algebra in
$mn$ variables $\C[\{x_{ij}\}], 1 \le i \le m, 1 \le j \le n$, the actions of
$\gl_m, \gl_n$ are given by
$$
E_{ab}^{\gl_m} =
\sum_{i=1}^n x_{ai}\frac\partial{\partial x_{bi}}, \qquad a,b = 1,\dots, m,
$$
$$
E_{ij}^{\gl_n} = \sum_{a=1}^m x_{ai} \frac \partial{\partial x_{aj}},
\qquad i,j = 1,\dots, n,
$$
where $E_{rs}^{\gl_k}$ are the standard generators of the Lie
algebra $\gl_k$, corresponding to the elementary matrices with $1$ at
the intersection of the $r$-th row and $s$-th column, and $0$ elsewhere.

We have the following theorem, see \cite{Zh, Howe}.

\begin{thm}
The $\gl_m \otimes\gl_n$ module $\mathbb P_{m,n}$ has the decomposition
$$
\mathbb P_{m,n} = \bigoplus_{\Lambda \in \mathcal P_{\min(m,n)}}
V_\Lambda^{\gl_m} \otimes V_\Lambda^{\gl_n}
$$
where $\mathcal P_j$ denotes the collection of finite sequences $\Lambda$ of nonnegative integers $\Lambda_1 \ge \dots \ge \Lambda_j$,
and $V_\Lambda^{\gl_k}$ for $k\ge j$ denotes the irreducible $\gl_k$-module with highest weight $(\Lambda_1,\dots,\Lambda_j,0,\dots,0)$.
\end{thm}

In particular, this theorem implies that a $\gl_n$-module $V_\Lambda^{\gl_n}$
may be realized as the subspace of $\gl_m$-singular vectors of weight $\Lambda$
in $\mathbb P_{m,n}$.

We will use the special case of the $(\gl_2,\gl_3)$ duality, which gives the following decomposition of the polynomial algebra $\mathbb P_{2,3}$:
$$
\mathbb P_{2,3} = \bigoplus_{m_1\ge m_2\ge 0} V_{(m_1,m_2)}^{\gl_2}
\otimes V_{(m_1,m_2,0)}^{\gl_3}.
$$

We identify a module $V_{(m_1,m_2,0)}^{\gl_3}$ with the the subspace of $\gl_2$-singular vectors in the algebra $\mathbb P_{2,3}$ of weight $(m_1,m_2)$ with respect to $\gl_2$:
$$V_{(m_1,m_2,0)}^{\gl_3} \cong \bigoplus_{l_1+l_2+l_3 = m_1+m_2} \( V_{(l_1,0)}^{\gl_2} \otimes V_{(l_2,0)}^{\gl_2} \otimes V_{(l_3,0)}^{\gl_2} \)^\sing[(m_1,m_2)].$$
Such identification gives the correspondence of the weight subspaces,
$$V_{(m_1,m_2,0)}^{\gl_3}[(l_1,l_2,l_3)] \cong \( V_{(l_1,0)}^{\gl_2} \otimes V_{(l_2,0)}^{\gl_2} \otimes V_{(l_3,0)}^{\gl_2} \)^\sing[(m_1,m_2)].$$

Our next goal is to describe a realization of $V_{(l_1,0)}^{\gl_2} \otimes
V_{(l_2,0)}^{\gl_2} \otimes V_{(l_3,0)}^{\gl_2}$ as a certain functional space.
In the remaining part of this section $\mathbf l = (l_1,l_2,l_3)$ will denote
a triple of non-negative integers.

To shorten formulas, we will omit superscripts for $\gl_2$-modules
and operators, and write $V_{(l,0)}, E_{ab}$ instead of $V_{(l,0)}^{\gl_2},
E_{ab}^{\gl_2}$, respectively. By a slight abuse of notation, we also write
$\1_l$ for the highest weight vector of the $\gl_2$-module $V_{(l,0)}$.

Let $\C[\mathbf z], \C(\mathbf z)$ denote respectively the algebras of
polynomial and rational functions in complex variables
$\mathbf z = (z_1,z_2,z_3)$.
%%%%%%%%%%%%%%%%%%%%%%%%%%%%%%%%%%%%%%%%%%%%%%%%%%%%%%%%%%%%%%%
%The symmetric group $S_3$ acts on these spaces
%by permutations of variables:
%$$
%(\sigma f)(z_1,z_2,z_3) = f(z_{\sigma(1)},z_{\sigma(2)},z_{\sigma(3)}),
%\qquad \sigma \in S_3\,.
%$$
%%%%%%%%%%%%%%%%%%%%%%%%%%%%%%%%%%%%%%%%%%%%%%%%%%%%%%%%%%%%%%%

Let $\mathcal H_k$ denote the space of symmetric polynomials in variables
$t_1,\dots,t_k$ of degree at most 2 in each $t_1,\dots,t_k$, with coefficients
in $\C[\mathbf z]$. We also set
$$
\mathcal H = \bigoplus_{k \in \Z_{\ge0}} \mathcal H_k\,.
$$

For any function $\varphi$ of $t_1,\dots,t_k$, we denote
$$
\Sym (\varphi (t_1,\dots,t_k) ) =
\sum_{\sigma \in S_k} \varphi(t_{\sigma(1)}, \dots, t_{\sigma(k)} )
$$

\begin{prop}\cite{TV3}
\label{thm:gl2H}
The space $\mathcal H$ has a structure of a $\gl_2$-module, depending on $\mathbf l$, such that the action of generators $\{ E_{ab} \}$ on a function $ \varphi \in \mathcal H_k$ is given by
\begin{gather*}
(E_{11} \varphi)(t_1,\dots,t_k) = ( l_1+l_2+l_3 - k ) \varphi(t_1,\dots,t_k), \label{eq:funcaction1}\\
(E_{22} \varphi)(t_1,\dots,t_k) = k \varphi(t_1,\dots,t_k),\\
(E_{12} \varphi)(t_1,\dots,t_{k-1}) = \lim_{t_k\to\infty} \frac{\varphi(t_1,\dots,t_k)}{t_k^2},
\\
\begin{aligned}
(E_{21} \varphi)(t_1,\dots,t_{k+1}) = \frac 1{k!}\;\Sym\( \varphi(t_1,\dots,t_k)
\( \prod_{a=1}^3 (t_{k+1}-z_a+l_a)
\prod_{i=1}^k \frac {t_{k+1}-t_i-1}{t_{k+1}-t_i} \,-{}\right.\right.\quad &
\\
\left.\left.\prod_{a=1}^3 (t_{k+1}-z_a)
\prod_{i=1}^k \frac {t_{k+1}-t_i+1}{t_{k+1}-t_i} \) \).\label{eq:funcaction4} &
\end{aligned}
\end{gather*}
\end{prop}

One can easily check that multiplication by any function of $\mathbf z$
commutes with the $\gl_2$ action above; thus, $\mathcal H$ can be regarded as a $\gl_2 \otimes \C[\mathbf z]$ - module.

\begin{prop}
\cite{TV3}
\label{thm:phimap}
The linear map
$$
\phi[\mathbf l]: V_{(l_1,0)} \otimes V_{(l_2,0)} \otimes V_{(l_3,0)} \to
\mathcal H\,,
$$
\begin{multline}
\phi[\mathbf l]\(E_{21}^{k_1} \1_{l_1}\otimes E_{21}^{k_2} \1_{l_2} \otimes
E_{21}^{k_3}\1_{l_3}\)\frac {l_1! l_2! l_3!}{(l_1-k_1)!(l_2-k_2)!(l_3-k_3)!}={}
\\
\Sym\(\prod_{a=1}^3 \(\prod_{i=1}^{k_1+\dots +k_{a-1}} (t_i-z_a)
\prod_{j=k_1+\dots+k_a+1}^{k_1+k_2+k_3} (t_j-z_a+l_a)\)
\prod_{i<j} \frac{t_i-t_j+1}{t_i-t_j}\),
\label{eq:phidef}
\end{multline}
is a homomorphism of $\gl_2$-modules.
\end{prop}

Let $\mathcal H_k[\mathbf l]$ denote the subspace of $\mathcal H_k$, consisting
of functions $\varphi \in \mathcal H_k$, satisfying additional "admissibility"
conditions
\begin{equation}
\varphi(z_i,z_i-1,\dots,z_i-l_i,t_{l_i+2},\dots,t_k) \equiv 0, \label{eq:adm}
\end{equation}
imposed if $k>l_i, \quad i =1,2,3$. We also denote
$$
\mathcal H[\mathbf l] = \bigoplus_{k\in\Z_{\ge0}} \mathcal H_k[\mathbf l]\,.
$$
One can check that $\mathcal H[\mathbf l]$ is
a $\gl_2 \otimes \C[\mathbf z]$-submodule of $\mathcal H$.

\begin{prop}\label{thm:phibasic}
\begin{enumerate}
\item
The map $\phi[\mathbf l]$ is injective, and
\begin{equation}
\operatorname{Im} \phi[\mathbf l] \otimes_{\C[\mathbf z]} \C(\mathbf z)=
\mathcal H[\mathbf l]\otimes_{\C[\mathbf z]} \C(\mathbf z).
\label{eq:dimensions}
\end{equation}
\item
There exists a homomorphism of $\gl_2 \otimes \C[\mathbf z]$-modules
$$
\mathcal I[\mathbf l]: \mathcal H[\mathbf l] \to
\(V_{(l_1,0)} \otimes V_{(l_2,0)} \otimes V_{(l_3,0)}\) \otimes \C(\mathbf z),
$$
such that
\begin{equation}
\(\mathcal I[\mathbf l] \circ \phi[\mathbf l]\) v = v,
\qquad\text{ for any }v\in V_{(l_1,0)} \otimes V_{(l_2,0)} \otimes V_{(l_3,0)}.
\label{eq:inversephi}
\end{equation}
Moreover, for any $\varphi \in \mathcal H[\mathbf l]$, the vector $\mathcal I[\mathbf l] \varphi$, rationally depending on $\mathbf z$, may have at most simple poles, located at the hyperplanes
$$z_j-r = z_i-l_i, \qquad r = 0,1,\dots,\min(l_i,l_j)-1, \quad 1 \le i \le j \le 3.$$
\end{enumerate}
\end{prop}
\noindent
The Proposition is proved in Appendix B.

\medskip
Introduce $\mathcal X[\mathbf l] \in \C[\mathbf z]$ by
$$
\mathcal X[\mathbf l; \mathbf z] =
\prod_{r=1}^{l_2} (z_1-l_1-z_2+r)
\prod_{r=1}^{l_3} ( z_1-l_1-z_3+r )
\prod_{r=1}^{l_3} ( z_2-l_2-z_3+r ),
$$
and set
\begin{equation}
\mathcal N[\mathbf l] = \mathcal X[\mathbf l] \ \mathcal I[\mathbf l].\label{eq:defN(z)}
\end{equation}

\begin{cor}\label{thm:polN(z)}
The formula \eqref{eq:defN(z)} defines a homomorphism of $\gl_2 \otimes \C[\mathbf z]$-modules
$$\mathcal N[\mathbf l]: \mathcal H[\mathbf l] \to \(V_{(l_1,0)} \otimes V_{(l_2,0)} \otimes V_{(l_3,0)}\) \otimes \C[\mathbf z].$$
\end{cor}

\begin{proof}
For any $\varphi \in \mathcal H[\mathbf l]$, the possible simple poles of
$\mathcal I[\mathbf l] \varphi$ are offset by $\mathcal X[\mathbf l]$.
Therefore, $\mathcal N[\mathbf l] \varphi$ is polynomial in $\mathbf z$.
The $\gl_2$-intertwining property for $\mathcal N[\mathbf l]$ follows from
the fact that $\mathcal I[\mathbf l]$ is a $\gl_2$-homomorphism, and
the commutativity of multiplication by any element from $\C[\mathbf z]$
with the action of $\gl_2$.
\end{proof}

Now we let formal variables $\mathbf z$ take particular complex values;
specifying $\mathbf z =\mathbf z_0$ for some $\mathbf z_0 \in \C^3$ will be
reflected by adding $\mathbf z_0$ to the notation. For example, $\phi[\mathbf
l; \mathbf z_0]$ will denote the composition of $\phi[\mathbf l]$ with the
homomorphism $\C[\mathbf z] \to \C$ of evaluation at $\mathbf z = \mathbf z_0$,
and $\mathcal H_k[\mathbf l; \mathbf z_0]$ will consist of symmetric
polynomials in $t_1,\ldots, t_k$ with complex coefficients, obtained from
elements of $\mathcal H_k[\mathbf l]$ by specializing $\mathbf z=\mathbf z_0$.

We now establish certain properties of operators $\phi[\mathbf l; \mathbf z_0], \mathcal N[\mathbf l; \mathbf z_0]$ for special values of $\mathbf z_0$.
Recall the following standard decompositions of the tensor product of $\gl_2$-modules:
$$V_{(l_1,0)} \otimes V_{(l_2,0)} \otimes V_{(l_3,0)} \cong
\(\bigoplus_{i=0}^{\min(l_1,l_2)} V_{(l_1+l_2-i,i)} \) \otimes V_{(l_3,0)} \cong
V_{(l_1,0)} \otimes \(\bigoplus_{i=0}^{\min(l_2,l_3)} V_{(l_2+l_3-i,i)} \).$$

\begin{prop}\label{thm:phiproperties}
Let $\mathbf z_0 \in \C^3$.
\begin{enumerate}
\item
If $\mathbf z_0$ is generic from the hyperplane $z_1-z_2-l_1+s = 0,\quad 0 \le s < \min(l_1,l_2)$,
then
$$
\ker \phi[\mathbf l; \mathbf z_0] =
\(\bigoplus_{i=s+1}^{\min(l_1,l_2)} V_{(l_1+l_2-i,i)} \) \otimes V_{(l_3,0)}\,.
$$
\item
If $\mathbf z_0$ is generic from the hyperplane $z_2-z_3-l_2+s = 0, \quad 0 \le s < \min(l_2,l_3)$,
then
$$
\ker \phi[\mathbf l; \mathbf z_0] =
V_{(l_1,0)} \otimes \(\bigoplus_{i=s+1}^{\min(l_2,l_3)} V_{(l_2+l_3-i,i)}\)\,,
$$
\end{enumerate}
\end{prop}

This Proposition will be proved in Appendix B. As a consequence, we obtain
the following properties of operators $\mathcal N[\mathbf l; \mathbf z_0]$.

\begin{cor}\label{thm:imageN(z)}
Let $\mathbf z_0 \in \C^3$.
\begin{enumerate}
\item
If $\mathbf z_0$ is generic from the hyperplane $z_1-l_1-z_2+s = 0, \quad 0 \le
s < \min(l_1,l_2)$, then
$$
\im \mathcal N[\mathbf l;\mathbf z_0] \subset
\(\bigoplus_{i=s+1}^{\min(l_1,l_2)} V_{(l_1+l_2-i,i)} \) \otimes V_{(l_3,0)}
$$
\item
If $\mathbf z_0$ is generic from the hyperplane $z_2-l_2-z_3+s = 0, \quad 0 \le
s < \min(l_2,l_3)$, then
$$
\im \mathcal N[\mathbf l;\mathbf z_0] \subset
V_{(l_1,0)} \otimes \(\bigoplus_{i=s+1}^{\min(l_2,l_3)} V_{(l_2+l_3-i,i)} \)
$$
\end{enumerate}
\end{cor}

\begin{proof}
Let $\varphi\in \mathcal H[\mathbf l; \mathbf z_0]$. Then it follows from the
definitions and \eqref{eq:dimensions} that
\begin{equation}
\phi[\mathbf l; \mathbf z_0] \circ \mathcal N[\mathbf l; \mathbf z_0] \varphi =
\mathcal X[\mathbf l; \mathbf z_0] \varphi\,.
\label{eq:kerim}
\end{equation}
For any point $\mathbf z_0$ from the hyperplane $z_1-l_1-z_2+s = 0$ with
$0 \le s < \min(l_1,l_2)$, we have that $\mathcal X[\mathbf l; \mathbf z_0]=0$,
and therefore $\mathcal N[\mathbf l;\mathbf z_0] \varphi \subset\ker\phi[\mathbf l; \mathbf z_0]$.
Hence, the first assertion follows from Proposition \ref{thm:phiproperties}.

The proof of the second assertion is similar.
\end{proof}

Note that admissibility conditions \eqref{eq:adm} are invariant under
simultaneous permutations of triples $\mathbf l$ and $\mathbf z$, i.e. we have
$$
\mathcal H[\mathbf l;\mathbf z] =
\mathcal H[\sigma(\mathbf l);\sigma(\mathbf z)], \qquad \sigma \in S_3.
$$

\begin{prop} \label{thm:Rmatrix}
There exists a family of linear operators
$\check {\mathcal R}_{l,l'}(u): V_{(l,0)} \otimes V_{(l',0)} \to
V_{(l',0)} \otimes V_{(l,0)}$, rationally depending on a complex parameter $u$,
such that
\begin{equation}
\check{\mathcal R}_{l,l'}(u) \1_{l} \otimes \1_{l'} =
\frac {\prod_{j=1}^{l'} (u-l+j) }{\prod_{j=1}^{l} (-u-l'+j)}
\1_{l'} \otimes \1_{l},
\label{eq:renormR}
\end{equation}
and the following diagram is commutative for generic $\mathbf z$:
\begin{equation}
\begin{CD}
\mathcal H[s_1(\mathbf l); s_1(\mathbf z)] @>{\mathcal N[s_1(\mathbf l);s_1(\mathbf z)]}>> V_{(l_2,0)} \otimes V_{(l_1,0)} \otimes V_{(l_3,0)} \\
@| @AA{\check{\mathcal R}_{l_1l_2}(z_1-z_2) \otimes 1}A\\
\mathcal H[\mathbf l; \mathbf z] @>{\mathcal N[\mathbf l; \mathbf z]}>> V_{(l_1,0)} \otimes V_{(l_2,0)} \otimes V_{(l_3,0)}
\\
@| @VV{1 \otimes \check{\mathcal R}_{l_2l_3}(z_2-z_3)}V\\
\mathcal H[s_2(\mathbf l);s_2(\mathbf z)] @>{\mathcal N[s_2(\mathbf l);s_2(\mathbf z)]}>> V_{(l_1,0)} \otimes V_{(l_3,0)} \otimes V_{(l_2,0)}
\end{CD}
\label{eq:N(z)}.
\end{equation}
\end{prop}

Proposition \ref{thm:Rmatrix} is proved in Appendix B, where it is shown that
the operators $\check{\mathcal R}_{l,l'}(u)$ are the standard rational
$R$-matrices, renormalized by \eqref{eq:renormR}.

%%%%

\section{The main construction}\label{sec:main}

Let $U$ be an irreducible finite-dimensional $\sl_3$-module with highest weight
$\Lambda$. Set
$$
m = \<\alpha_1+2\alpha_2,\Lambda\>, \qquad k = \<\alpha_2,\Lambda\>.
$$
The $\gl_3$-module $V_{(m-k,k,0)}^{\gl_3}$, regarded by restriction as an
$\sl_3$-module, is isomorphic to $U$. We fix an identification
(which is unique up to a scalar)
$$
\Theta: U \to \bigoplus_{\mathbf l \in \Z_{\ge0}^3}
\(V_{(l_1,0)} \otimes V_{(l_2,0)} \otimes V_{(l_3,0)}\)^\sing[(m-k,k)],
$$
so that $\Theta$ maps a weight subspace $U[\mu]$ to $\(V_{(l_1,0)} \otimes
V_{(l_2,0)} \otimes V_{(l_3,0)}\)^\sing[(m-k,k)]$, where the corresponding $\mathbf l =
(l_1,l_2,l_3)\in\Z_{\ge0}^3$ is determined by
\begin{equation}
l_1-l_2 = \<\alpha_1,\mu\>, \qquad l_2-l_3 = \<\alpha_2,\mu\>, \qquad l_1+l_2+l_3 = m. \label{eq:lmu}
\end{equation}
Under this identification the Weyl group for $\sl_3$ coincides with the symmetric
group $S_3$ naturally acting on $\Z_{\ge0}^3$.

For any $\mathbf l\in\Z_{\ge0}^3$ we fix a basis $\{v_i[\mathbf l] \}$ of the
space $\(V_{(l_1,0)}\otimes V_{(l_2,0)}\otimes V_{(l_3,0)}\)^\sing[(m-k,k)]$.
We assume that $v_i[w(\mathbf l)] = \bar w v_i[\mathbf l]$, where the operators
$\bar w$ correspond to the action of the Weyl group given by Lemma \ref{thm:wbarep}.
Note that if $w(\mathbf l)=\mathbf l$, then the restriction of $\bar w$
on $\(V_{(l_1,0)}\otimes V_{(l_2,0)}\otimes V_{(l_3,0)}\)^\sing[(m-k,k)]$
is the identity operator.

We have the following Lemma, proved in Appendix C.

\begin{lem}\label{thm:constantbasis}
For any $\mathbf l = (l_1,l_2,l_3) \in \Z_{\ge0}^3$, and $k \in \Z_{\ge0}$, we have
$$
\dim \mathcal H_k[\mathbf l; \mathbf z]^\sing =
\dim \(V_{(l_1,0)} \otimes V_{(l_2,0)} \otimes V_{(l_3,0)}\)^\sing[(m-k,k)]\,,
$$
and there exists a $\C[\mathbf z]$-basis $\{ \varsigma_i[\mathbf l] \}$ of
$\mathcal H_k[\mathbf l]^\sing$, such that
\begin{equation}
\varsigma_i[\mathbf l; \mathbf z] = \varsigma_i[\sigma(\mathbf l);\sigma(\mathbf z)], \qquad \sigma \in S_3.
\label{eq:invaradm}
\end{equation}
\end{lem}

Define the operators
$\Upsilon[\mathbf l]: \(V_{(l_1,0)} \otimes V_{(l_2,0)} \otimes V_{(l_3,0)}\)^\sing[(m-k,k)] \to \mathcal H[\mathbf l]^\sing$ by
$$\Upsilon[\mathbf l] v_i[\mathbf l] = \varsigma_i[\mathbf l].$$
By construction, the operators $\Upsilon[\mathbf l]$ satisfy
\begin{equation}
\begin{CD}
 \(V_{(l_2,0)} \otimes V_{(l_1,0)} \otimes V_{(l_3,0)}\)^\sing[(m-k,k)] @>{\Upsilon[s_1(\mathbf l); s_1(\mathbf z)]}>> \mathcal H_k[s_1(\mathbf l); s_1(\mathbf z)]^\sing \\
@A{\bar s_1}AA @|\\
 \(V_{(l_1,0)} \otimes V_{(l_2,0)} \otimes V_{(l_3,0)}\)^\sing[(m-k,k)] @>{\Upsilon[\mathbf l; \mathbf z]}>> \mathcal H_k[\mathbf l; \mathbf z]^\sing \\
@V{\bar s_2}VV @|\\
 \(V_{(l_1,0)} \otimes V_{(l_3,0)} \otimes V_{(l_2,0)}\)^\sing[(m-k,k)] @>{\Upsilon[s_2(\mathbf l); s_2(\mathbf z)]}>> \mathcal H_k[s_2(\mathbf l);s_2(\mathbf z)]^\sing
\end{CD}
\label{eq:Upsilon}.
\end{equation}

We now introduce our main object --- the operators
$N_\mu(\lambda)\in\End_\C(U[\mu])$ --- by setting
\begin{equation}
N_\mu(\lambda) = \ \Theta^{-1} \circ\mathcal N[\mathbf l;\mathbf z] \circ
\Upsilon[\mathbf l;\mathbf z] \circ\Theta\,.
\end{equation}
Here $\mu$ and $l_1,l_2,l_3$ are as in \eqref{eq:lmu},
and $z_1,z_2,z_3$ satisfy
\begin{equation}
(z_1-l_1)-(z_2-l_2) = \<\alpha_1,\lambda+\rho\>, \qquad
(z_2-l_2)-(z_3-l_3) = \<\alpha_2,\lambda+\rho\>\,.
\label{eq:zlambda}
\end{equation}

\begin{thm}
The operator $N(\lambda) = \bigoplus_\mu N_\mu(\lambda)$ is a regularizing
operator for the $\g$-module $U$. The corresponding Weyl group representation
$w \mapsto a_w$ is such that $a_{s_i} = \bar s_i$.
\end{thm}

\begin{proof}
Corollary \ref{thm:polN(z)} implies that the operator $N(\lambda)$ depends
on $\lambda \in \h^*$ polynomially.

To check \eqref{eq:dynN}, we use the result of \cite{TV2}, where it was shown
that the action of the dynamic Weyl group operators $A_w$ in $\mathbb P_{2,3}$
coincides with the action of operators $\check {\mathcal R}$. More precisely,
we have the following commutative diagram:
\begin{equation}
\begin{CD}
U[s_1\mu] @>>> \(V_{(l_2,0)} \otimes V_{(l_1,0)} \otimes V_{(l_3,0)}\)^\sing[(m-k,k)]\\
@A{A_{s_1}^\mu(\lambda)}AA @AA{\check {\mathcal R}_{l_1l_2}(z_1-z_2) \otimes 1}A\\
U[\mu] @>>> \(V_{(l_1,0)} \otimes V_{(l_2,0)} \otimes V_{(l_3,0)}\)^\sing[(m-k,k)]\\
@V{A_{s_2}^\mu(\lambda)}VV @VV{1 \otimes \check{\mathcal R}_{l_2l_3}(z_2-z_3)}V\\
U[s_2\mu] @>>> \(V_{(l_1,0)} \otimes V_{(l_3,0)} \otimes V_{(l_2,0)}\)^\sing[(m-k,k)]
\end{CD}\ .
\label{eq:dynRmatrix}
\end{equation}
Combining \eqref{eq:dynRmatrix} with \eqref{eq:N(z)} and \eqref{eq:Upsilon},
we get the commutative diagram
$$
\begin{CD}
U[s_1 \mu] @>{N_{s_1\mu}(s_1 \cdot \lambda)}>> U[s_1 \mu] \\
@A{\bar s_1}AA @AA{A_{s_1}^\mu(\lambda)}A\\
U[\mu] @>{N_\mu(\lambda)}>> U[\mu] \\
@V{\bar s_2}VV @VV{A_{s_2}^\mu(\lambda)}V\\
U[s_2 \mu] @>{N_{s_2 \mu}(s_2\cdot\lambda)}>> U[s_2 \mu]
\end{CD}\quad\ ,
$$
which implies \eqref{eq:dynN}.

Next, the operator $\mathcal N[\mathbf l; \mathbf z]^{-1} =
\mathcal X[\mathbf l; \mathbf z]^{-1}\phi[\mathbf l; \mathbf z]$
is singular only at the hyperplanes
$$
z_1-l_1-z_2+s = 0, \qquad 0 \le s < l_2\,,
$$
$$
z_2-l_2-z_3+s = 0, \qquad 0 \le s < l_3\,,
$$
$$
z_1-l_1-z_3+s = 0, \qquad 0 \le s < l_3\,,
$$
which correspond respectively to the hyperplanes
$$
\chi_{\alpha_1,l_2-s}(\lambda) = 0\,,\qquad 0 \le s < l_2\,,
$$
$$
\chi_{\alpha_2,l_3-s}(\lambda) = 0\,,\qquad 0 \le s < l_3\,,
$$
$$
\chi_{\alpha_1+\alpha_2,l_3-s}(\lambda) = 0\qquad 0 \le s < l_3\,.
$$
Now it follows from Lemma \ref{thm:detfrompoles} that $N(\lambda)$ satisfies
\eqref{eq:detN}.

Finally, let us show that for any $u \in U[\mu]$, the vector $\Xi_N(\lambda)u
\in \mathcal U(\n^-) \otimes U$ polynomially depends on $\lambda$. According
to \eqref{eq:Xi}, this is equivalent to regularity of the sum
$\sum_j \(S_\lambda^{-1}\)_{ij} \omega(g_j) N(\lambda)u$ for any $i$ at any
hyperplane $\chi_{\alpha,r}(\lambda) = 0,\quad \alpha\in\Delta^+, r = 1,2, \dots$.

Consider the case $\alpha = \alpha_1$.

Let $\{g_j\}$ be the Poincare-Birkhoff-Witt basis of $\mathcal U(\n^-)$,
corresponding to an ordering of $\Delta^+$ with $\alpha_1$ being the last root.
It is known (see \cite{ES}) that $\(S_\lambda^{-1}\)_{ij}$ is regular at the
hyperplane $\chi_{\alpha_1,r}(\lambda)=0$ unless $g_j = \tilde g_j f_1^r$ for
some $\tilde g_j \in \mathcal U(\n^-)$. In the last case
$\(S_\lambda^{-1}\)_{ij}$ may have a simple pole there, and the desired
regularity will follow if $\omega(\tilde g_j) e_1^r N(\lambda)u$ vanishes at
the hyperplane $\chi_{\alpha_1,r}(\lambda)=0$.

We prove that $e_1^r N(\lambda)u=0$ at the hyperplane
$\chi_{\alpha_1,r}(\lambda) = 0$ using the explicit realization $\Theta$
of the module $U$ in $\mathbb P_{2,3}$. The operator $e_1^r: U[\mu] \to
U[\mu+r\alpha_1]$ corresponds to the operator
$$
\(E_{12}^{\gl_3}\)^r:
\( V_{(l_1,0)} \otimes V_{(l_2,0)} \otimes V_{(l_3,0)} \)^\sing[(m-k,k)] \to
\( V_{(l_1+r,0)} \otimes V_{(l_2-r,0)} \otimes V_{(l_3,0)} \)^\sing[(m-k,k)]\,,
$$
if $r\le l_2$, and is equal to zero, if $r>l_2$. There is nothing to prove in
the latter case. The claim in the former case is implied by the following
lemma.

\begin{lem}
Assume that $r\le l_2$. Let $\mathbf z_0$ be a generic point from the
hyperplane $(z_1-l_1)-z_2+(l_2-r) = 0$. Then the following composition
of operators vanishes:
\begin{equation}
\mathcal H[\mathbf l; \mathbf z_0] \stackrel{\mathcal N[\mathbf l; \mathbf z_0]}{\longrightarrow}
V_{(l_1,0)} \otimes V_{(l_2,0)} \otimes V_{(l_3,0)} \stackrel{\(E_{12}^{\gl_3}\)^{r} }{\longrightarrow}
V_{(l_1+r,0)} \otimes V_{(l_2-r,0)} \otimes V_{(l_3,0)}.
\label{eq:composition}
\end{equation}
\end{lem}

\begin{proof}
Since $r\le l_2$, the product $\mathcal X[\mathbf l; \mathbf z]$ contains
the factor $z_1-l_1-z_2+(l_2-r)$, corresponding to
$\chi_{\alpha_1,r}(\lambda)$, and $\mathcal X[\mathbf l; \mathbf z_0] = 0$
If $l_2-r > \min(l_1,l_2)$, then $\mathcal I[\mathbf l]\varphi$ is regular
at $\mathbf z = \mathbf z_0$ for any $\varphi\in\mathcal H[\mathbf l]$, and
the operator $\mathcal N[\mathbf l; \mathbf z_0] =
\mathcal X[\mathbf l;\mathbf z_0]\,\mathcal I[\mathbf l;\mathbf z_0]$ vanishes.

If $l_2-r \le \min(l_1,l_2)$, then by Proposition \ref{thm:imageN(z)},
the image of $\mathcal N[\mathbf l; \mathbf z_0]$ is contained in $\(
\bigoplus_{i=l_2-r+1}^{\min(l_1,l_2)}V_{(l_1+l_2-i,i)}\)\otimes V_{(l_3,0)}$.
Since $\(E_{12}^{\gl_3}\)^{r}$ commutes with the $\gl_2$ action and
$$
V_{(l_1+r,0)} \otimes V_{(l_2-r,0)} \otimes V_{(l_3,0)} \cong
\(\bigoplus_{i=0}^{l_2-r} V_{(l_1+l_2-i,i)} \) \otimes V_{(l_3,0)}\,,
$$
the composition \eqref{eq:composition} vanishes.
\end{proof}

The case $\alpha = \alpha_2$ can be considered similarly to the case
$\alpha = \alpha_1$.

Let now $\alpha = \alpha_1 + \alpha_2$. In this case we have $\alpha = s_2 \alpha_1$, and therefore
$$
\chi_{\alpha,r}(\lambda) = \<\a, \lambda+\rho\> - r =
\<s_2 \a_1, \lambda+\rho\> - r = \<\a_1,s_2 \cdot \lambda + \rho \> - r = \chi_{\alpha_1,r}(s_2 \cdot \lambda).
$$

Recall that $\Xi_N(\lambda)u$ is regular at the hyperplane
$\chi_{\a,r}(\lambda)=0$ if and only if the operator
$\Phi_\lambda^{N(\lambda)u}$ is regular at the same hyperplane.

Consider $\lambda$ from the hyperplane $\<\a_2, \lambda+\rho\>=-n$, for
a large enough nonnegative integer $n$. Then $M_\lambda$ is contained as
a proper submodule in $M_{s_2\cdot\lambda}$, and we have
$$
\Phi_\lambda^{N(\lambda)u} =
\Phi_{s_2\cdot\lambda}^{A(s_2\cdot\lambda)N(s_2\cdot\lambda)u}|_{M_\lambda}\,.
$$
Since the operator $\Phi_{s_2\cdot\lambda}^{N(s_i\cdot\lambda)u}$ is regular
at the hyperplanes $\chi_{\a_1,r}(s_2\cdot\lambda)=0$, so is the operator
$\Phi_\lambda^{N(\lambda)u}$. Hence $\Xi_N(\lambda)u$ has no pole
at the hyperplane $\chi_{\a,r}(\lambda)=\chi_{\a_1,r}(s_2\cdot\lambda)=0$,
provided $\<\a_2, \lambda+\rho\>=-n$.

The rational function $\chi_{\a,r}(\lambda)\Xi_N(\lambda)u$ is regular at
the hyperplane $\chi_{\a,r}(\lambda)=0$, and vanishes at infinitely many points
of the hyperplane. Therefore, the function must vanish at the hyperplane
identically, that is $\Xi_N(\lambda)u$ has no pole at the hyperplane
$\chi_{\a,r}(\lambda)=0$.
\end{proof}

%%%%

\section{Resonance conditions}

Let $\g$ be a semisimple Lie algebra, and let $U$ be an irreducible
finite-dimensional $\g$-module with nontrivial zero weight subspace $U[0]$.

Let $N(\lambda) \in \End_\C(U)$ be a regularizing operator.
Introduce an $\End_\C(U[0])$-valued function $\Psi(\lambda,x)$, where $\lambda \in \h^*, x \in \h$, by
$$\Psi(\lambda,x) u = \operatorname{Tr}|_{M_\lambda} \( \Phi_\lambda^{N(\lambda)u} e^x \).$$

Special cases of the function $\Psi(\lambda,x)$ were used in \cite{ES} to give
a representation-theoretic proof of the algebraic integrability of
Calogero-Sutherland systems. In \cite{EV} it was shown that the trace functions
satisfy remarkable differential and difference equations.

An important property of the function $\Psi(\lambda,x)$ is that it satisfies
certain resonance conditions.

\begin{prop}
Let $N = N(\lambda)$ be a regularizing operator acting in a finite-dimensional
$\g$-module $U$ with $U[0]\ne 0$, and let $\Psi(\lambda,x)$ be the
corresponding trace function.

Let $\a\in\Delta^+$ and $k\in\Z_{\ge0}$ be such that $U[k\alpha]\ne 0$.
Then for any $\lambda \in \h^*$ such that $\chi_{\alpha,k}(\lambda)=0$,
and any $u\in\ker N(\lambda) \subset U[0]$, we have
$$
\Psi (\lambda,x) u = \Psi(s_\alpha\cdot\lambda,x)u\,.
$$
\end{prop}

\begin{proof}
Let $\lambda$ be generic from the hyperplane $\chi_{\alpha,k}(\lambda)=0$. Then $M_\lambda$ contains a unique proper submodule $M_{s_\alpha \cdot \lambda}$, and $N(\lambda)u = 0$ implies that the image of
$\Phi_\lambda^{N(\lambda)u}$ is contained in $M_{s_\alpha \cdot \lambda} \otimes U$. Therefore,
the trace of the operator $\Phi_\lambda^{N(\lambda)u}$ in $M_\lambda$ is equal to the trace of its resrtiction to $M_{s_\alpha \cdot \lambda}$.
Moreover, the restriction of $\Phi_\lambda^{N(\lambda)u}$ to the submodule $M_{s_\alpha \cdot \lambda}$
coincides with $\Phi_{s_\alpha \cdot \lambda}^{N(s_\alpha \cdot\lambda)u}$, and we have
\begin{align*}
\Psi(\lambda,x)u =
\operatorname{Tr}|_{M_\lambda} \( \Phi_\lambda^{N(\lambda)u} e^x \)=
\operatorname{Tr}|_{M_{s_\alpha\cdot\lambda}}
\( \Phi_\lambda^{N(\lambda)u} e^x \) ={}&
\\[3pt]
\operatorname{Tr}|_{M_{s_\alpha\cdot\lambda}}
\( \Phi_{s_\alpha\cdot\lambda}^{N(s_\alpha\cdot\lambda)u} e^x \) =
\Psi(s_\alpha \cdot \lambda,x) u\,. &
\end{align*}
\end{proof}

{\bf Example:}\enspace
Let $\g=\sl_2$. We identify $\h^*$ with $\C$ by associating $ \h^* \ni \lambda
\leftrightarrow \<\alpha,\lambda\> \in \C$.

Let $U$ be an irreducible finite-dimensional $\sl_2$-module with even highest
weight $\Lambda$. In this case the weight subspace $U[0]$ is one-dimensional,
$\End(U[0])\cong\C$, and the function $\Psi(\lambda,x)$ is scalar-valued.
It satisfies the resonance conditions
$$
\Psi(\lambda,x) = \Psi(s\cdot\lambda,x), \qquad
\lambda = 0,1,2,\dots, \Lambda/2-1\,.
$$

%%%%

\appendix

\section{Formal monomials}

To prove Theorem \ref{thm:dyngroup} we use the calculus of formal monomials in $\mathcal U(\n^-)$, developed in \cite{FFM}.

\begin{lem}
Let $f \in \n^-, X \in \g$. For any $n \in \Z_{\ge0}$ we have the following identity in $U(\n^-)$,
\begin{equation}
X f^n = \sum_{k=0}^\infty \binom{n}{k} f^{n-k} \underbrace{[[[X,f],\dots,],f]}_{k \text { times }}.
\label{eq:commutation}
\end{equation}
\end{lem}
This lemma is proved by induction. Note that the summation over $k$ is finite.

Introduce an associative algebra $\mathbb A$, with unit $\1$ and generators
$f_i^\gamma$, $i = 1,\dots,\dim\h$, $\gamma \in \C$, subject to the
relations
\begin{gather}
\underbrace{[f_i,[,\dots,[f_i,f_j]\dots]]}_{a_{ij} \text { times }} = 0, \qquad
a_{ij} = \frac{2\<\alpha_i,\alpha_j\>}{\<\alpha_i,\alpha_i\>},\label{eq:Serre} \\
f_i^{\beta} f_i^\gamma = f_i^{\beta+\gamma}, \qquad \beta,\gamma \in \C,\\
f_i f_j^\gamma = \sum_{k=0}^\infty \binom{\gamma}{k} f_j^{\gamma-k}
\underbrace{[[\dots[f_i,f_j],\dots,],f_j]}_{k \text { times }}, \qquad \gamma \in \C.
\label{eq:formal}
\end{gather}
where the commutator $[\cdot,\cdot]$ is defined as $[X,Y] = XY - YX$.

\begin{lem}
There exists an algebra homomorphism $\iota: \mathcal U(\n^-) \to \mathbb A$ such that
$\iota(f_i) = f_i$ for any $i=1,\dots,\dim\h$. The homomorphism $\iota$ is injective.
\end{lem}
\begin{proof}
It is known that $U(\n^-)$ is the quotient of the free Lie algebra with
generators $\{f_i\}$, by the ideal, generated by the Serre relations
\eqref{eq:Serre}. Therefore, $\iota$ extends from the generators to the entire
algebra $\mathcal U(\n^-)$.

Next, we note that the algebra $\mathbb A$ is nonzero, i.e. the ideal,
generated by the relations \eqref{eq:formal} is proper. Indeed, a linear
functional $\varphi$ on $\mathbb A$, defined on monomials by
$$
\varphi(f_{i_1}^{\gamma_1} \dots f_{i_l}^{\gamma_l}) = 1,
$$
is well-defined, and is nonzero, so $\mathbb A \ne 0$.

Let $\lambda \in \h^*$ be generic. Introduce a $\g$-module $\tilde M_\lambda$,
which is equal to $\mathbb A$ as a vector space, with $f_i$ acting naturally on
the left, $h \in \h$ acting by
$$
h f_{i_1}^{\gamma_1} \dots f_{i_l}^{\gamma_l} =
\bigl\<h,\lambda - \sum_{k=1}^l \gamma_k \alpha_{i_k} \bigr\>
f_{i_1}^{\gamma_1} \dots f_{i_l}^{\gamma_l}.
$$
The action of operators $e_i$ is determined by the condition $e_i \1 = 0$, and the commutation formula
$$
e_i f_j^\gamma = \sum_{k=0}^\infty \binom{\gamma}{k} f_j^{\gamma-k}
\underbrace{[[\dots[e_i,f_j],\dots,],f_j]}_{k \text { times }},
\qquad \gamma \in \C.
$$
The module $\tilde M_\lambda$ is an obvious analog of the Verma module
$M_\lambda$.

For any $X,Y \in \mathcal U(\n^-)$, consider the action of $\omega (Y) \in
\mathcal U(\g)$ on the element $\iota(X)\1 \in \tilde M_\lambda$. We have
$\omega(Y) \iota(X) \1 = S_\lambda(X,Y) \1$, and therefore $\iota(X) = 0$
implies $S_\lambda(X,Y) = 0$ for any $Y \in \mathcal U(\n^-)$. Since
the Shapovalov form $S_\lambda$ is nondegenerate for generic $\lambda$,
this means that $X = 0$. Thus, $\iota$ is an injection.
\end{proof}

We identify $\mathcal U(\n^-)$ with its image $\iota(\mathcal U(\n^-)) \subset \mathbb A$.
We say that an element $X \in \mathbb A$ makes sense, if $X \in \mathcal U(\n^-)$.

\begin{example}
Let $\g = \sl_3$, and let $\gamma_1,\gamma_2,\gamma_3 \in \C$. If
$\gamma_1+\gamma_3, \gamma_2\in\Z_{\ge0}$, and $\gamma_1+\gamma_3\ge\gamma_2$,
then the monomial $f_1^{\gamma_1} f_2^{\gamma_2} f_1^{\gamma_3}$ makes sense,
and we can write
$$
f_1^{\gamma_1} f_2^{\gamma_2} f_1^{\gamma_3} =
\sum_{j=0}^\infty \binom{\gamma_2}{j} \binom{\gamma_3}{j}
f_1^{\gamma_1+\gamma_3-j} f_2^{\gamma_2-j} [f_2,f_1]^j.
$$
\end{example}

Consider the completion $\mathbb A \hat{{}\otimes{}} \mathcal U(\n^-)$, consisting of
possibly infinite sums $\sum_i X_i \otimes Y_i$ with $X_i \in \mathbb A$ and
homogeneous $Y_i \in \mathcal U(\n^-)$, such that $\wt (Y_i) \to \infty$ when
$i \to \infty$.

We following two Lemmas can be easily verified; we give sketches of proof,
and leave more technical details to a reader.

\begin{lem}
There exists an algebra homomorphism
$$\Delta: \mathbb A \to \mathbb A \hat{{}\otimes{}} \mathcal U(\n^-),$$
such that
$$
\Delta(f_i^\gamma) = \sum_{j=0}^\infty \binom {\gamma}{j}
f_i^{\gamma-j} \otimes f_{i}^{j}.
$$
\qed
\end{lem}
\begin{proof}[Sketch of proof]
We have to check that the comultiplication above can be extended to an algebra
homomorphism, i.e.~that it respects the defining relations \eqref{eq:formal}.
For each such relation of homogeneous weight $\gamma$, it amounts to
verifying a family of identities in weight subspaces $\mathbb A[\nu] \otimes
\mathcal U(\n^-)[\gamma-\nu]$. These identities are are valid for any nonnegative
integer $\gamma$, because in this case they are equivalent to the consistency of
the comultiplication in $\mathcal U(\n^-)$. Since these identities are
polynomial in $\gamma$, they are satisfied for arbitrary $\gamma$.
\end{proof}

It is easy to see that if $X \in \mathcal U(\n^-) \subset \mathbb A$, then
$\Delta(X) \in \mathcal U(\n^-) \otimes \mathcal U(\n^-) \subset \mathbb A
\hat{{}\otimes{}} \mathcal U(\n^-)$.

\begin{lem}\label{lem:PBW}
Let $J \in \Z_{\ge0}, \gamma_1,\dots,\gamma_l \in \C$. Let $\{g_i\}$ be
a homogeneous basis of $\mathcal U(\n^-)$. Then there exist
$J_1,\dots,J_l \in \Z_{\ge0}$, such that any monomial $f_{i_1}^{\gamma_1-j_1}
\dots f_{i_l}^{\gamma_l-j_l}$ with $j_1,\dots,j_l \le J$, can be written as
$$
f_{i_1}^{\gamma_1-j_1} \dots f_{i_l}^{\gamma_l-j_l} =
f_{i_1}^{\gamma_1-J_1} \dots f_{i_l}^{\gamma_l-J_l}
\sum_i q_i(\gamma_1,\dots,\gamma_l) g_i,
$$
for some polynomials $q_i(\gamma_1,\dots,\gamma_l)$.
\qed
\end{lem}

\begin{proof}[Sketch of proof]
The statement is obvious for $l=1$. Using induction on $l$, we can write
$$
f_{i_1}^{\gamma_1-j_1}\dots f_{i_{l-1}}^{\gamma_{l-1}-j_{l-1}}
f_{i_l}^{\gamma_l-j_l} =
f_{i_1}^{\gamma_1-J_1} \dots f_{i_{l-1}}^{\gamma_{l-1}-J_{l-1}}
\( \sum_i q_i(\gamma_1,\dots,\gamma_l) g_i \) f_{i_l}^{\gamma_l-j_l}.
$$
Then we can move all $f_{i_l}$ to the left by commuting them with $\{g_i\}$;
nilpotency of the adjoint action of $f_{i_l}$ will guarantee that all only
finite number of $f_{i_l}$'s can be absorbed.
\end{proof}

We now use these resulst to give a proof of Theorem \ref{thm:dyngroup}.

\begin{proof}[Proof of Theorem \ref{thm:dyngroup}]

We represent the left side of \eqref{eq:dyngroup} as
\begin{align}
& \Phi_\lambda^u v_{w\cdot(\lambda+\mu)} =
\Phi_\lambda^u F_w(\lambda+\mu) \1_{\lambda+\mu} =
\Delta(F_w(\lambda+\mu))\Xi(\lambda) (\1_\lambda \otimes u) =
\notag
\\[4pt]
& \sum_{j_1,\dots,j_l \in \Z_{\ge0}}\sum_{i,j}
q_{j_1,\ldots,j_l}(\gamma_1(\lambda+\mu),\ldots,\gamma_l(\lambda+\mu))\times{}
\notag
\\
& \qquad f_{i_1}^{\gamma_1(\lambda+\mu)-j_1} \dots
f_{i_l}^{\gamma_l(\lambda+\mu)-j_l} g_i \1_\lambda \otimes
\(S_\lambda^{-1}\)_{ij} f_{i_1}^{j_1} \dots f_{i_l}^{j_l} \omega(g_j) u.
\label{eq:DeltaF}
\end{align}
Here the coefficients
$q_{j_1,\ldots,j_l}(\gamma_1(\lambda+\mu),\ldots,\gamma_l(\lambda+\mu))$ are
certain polynomials in $\gamma_1(\lambda+\mu),\ldots,\gamma_l(\lambda+\mu)$.
The summation in \eqref{eq:DeltaF} is actually finite, because the element
$f_{i_1}^{j_1} \dots f_{i_l}^{j_l}$ acts as zero in $U$ if any of
$j_1,\dots,j_l$ is big enough.

Next, we have that $\gamma_k(\lambda+\mu)=\gamma_k(\lambda) + \<\beta_k,\mu\>$
for some $\beta_k \in \Delta^+$. Hence, for all nonzero monomials in the sum,
we have
$$
\gamma_k(\lambda+\mu) - j_k =
\gamma_k(\lambda) + \<\beta_k,\mu\> - j_k \ge \gamma_k(\lambda)-J,
\qquad k = 1,\dots,l,
$$
for some $J \in \Z_{\ge0}$. Lemma \ref{lem:PBW} implies that there exist $J_1,\dots,J_l \in \Z_{\ge 0}$, such that we can write
\begin{equation}
\Delta(F_w(\lambda+\mu))\Xi(\lambda) = \(f_{i_1}^{\gamma_1(\lambda)-J_1} \dots f_{i_l}^{\gamma_l(\lambda)-J_l} \otimes 1 \)
\(\sum_i g_i \otimes P_i(\lambda)\),
\label{eq:comult}
\end{equation}
for a certain finite collection of $\End_\C(U)$-valued rational functions
$P_i(\lambda)$.

We now make two claims, which imply Theorem \ref{thm:dyngroup}.

First, in the summation formula \eqref{eq:comult} we have
\begin{equation}
P_i(\lambda) \equiv 0, \qquad \text{ if \ }
\wt(g_i) \notin \wt(f_{i_1}^{J_1} \dots f_{i_l}^{J_l}) - Q^+.
\label{eq:vanish}
\end{equation}

Second, for some rational $\End_\C(U)$-valued functions $A_w(\lambda)$, we have
\begin{equation}
\(f_{i_1}^{\gamma_1(\lambda)-J_1} \dots f_{i_l}^{\gamma_l(\lambda)-J_l} \otimes 1 \)
\( {\sum_i}' g_i \otimes P_i(\lambda)\) = F_w(\lambda) \otimes A_w(\lambda),
\label{eq:claim2}
\end{equation}
where ${\sum_i}'$ denotes the sum of terms with $\wt(g_i) = \wt(f_{i_1}^{J_1} \dots f_{i_l}^{J_l})$.

To verify these claims, we consider a dominant integral weight $\lambda$ such
that $\lambda \gg0$. One can check that for any $u\in U$ we have
\begin{equation}
\Delta(F_w(\lambda+\mu))\Xi(\lambda) (\1_\lambda \otimes u)
= v_{w\cdot\lambda} \otimes \tilde u + \text{ lower order terms },
\label{eq:biglambda}
\end{equation}
for some $\tilde u \in U$.
In particular, the term $g_i \otimes P_i(\lambda)$ is nonzero only if
\begin{equation}
\wt (f_{i_1}^{\gamma_1(\lambda)-J_1} \dots f_{i_l}^{\gamma_l(\lambda)-J_l}) +
\wt(g_i) + \lambda \le w\cdot \lambda.
\end{equation}
{}From the definition of $\gamma_k(\lambda)$ it follows that
$$\lambda + \wt (f_{i_1}^{\gamma_1(\lambda)} \dots f_{i_l}^{\gamma_l(\lambda)}) = w\cdot \lambda,$$
and thus \eqref{eq:biglambda} implies \eqref{eq:vanish} for the given $\lambda$.
Since the rational functions $P_i(\lambda)$ vanish for all sufficiently large
dominant integral weights, they must vanish for all $\lambda$.

Using equation \eqref{eq:claim2}, we define the dynamical Weyl group operator
$A_w(\lambda)$ for integral dominant $\lambda \gg0$ by the rule
$A_w(\lambda) u = \tilde u$. We now establish that $A_w(\lambda)$ rationally
depends on $\lambda$, and verify our second claim.

By Lemma \ref{lem:PBW}, we can write
$$F_w(\lambda) f_{i_1}^{\gamma_1(\lambda)-J_1} \dots f_{i_l}^{\gamma_l(\lambda)-J_l} \sum_i Q_i(\lambda) g_i,$$
for some polynomials $Q_i(\lambda)$ in $\lambda \in \h^*$. Equation \eqref{eq:biglambda} then becomes
\begin{multline}
\(f_{i_1}^{\gamma_1(\lambda)-J_1} \dots f_{i_l}^{\gamma_l(\lambda)-J_l} \otimes 1\)
\( {\sum_i}' g_i \1 \otimes P_i(\lambda)u \) = \\
\(f_{i_1}^{\gamma_1(\lambda)-J_1} \dots f_{i_l}^{\gamma_l(\lambda)-J_l} \otimes 1\)
\( {\sum_i}' Q_i(\lambda) g_i \1 \otimes A_w(\lambda) u\).
\label{eq:A(w)}
\end{multline}
or equivalently
$$ {\sum_i}' g_i \1 \otimes P_i(\lambda)u = {\sum_i}' Q_i(\lambda) g_i \1 \otimes A_w(\lambda) u.$$

It is valid for dominant integral $\lambda\gg0$, and implies that for all $i$ we have
\begin{equation}
A_w(\lambda) = \frac {P_i(\lambda)}{Q_i(\lambda)}, \qquad \text{ if }
\wt(g_i) = \wt(f_{i_1}^{J_1} \dots f_{i_l}^{J_l}).
\label{eq:Awratio}
\end{equation}
Therefore, $A_w(\lambda)$ can be extended to a rational $\End_\C(U)$-valued
function, such that \eqref{eq:Awratio} is satisfied for all $\lambda$.
Equations \eqref{eq:A(w)}, \eqref{eq:claim2}, and \eqref{eq:dyngroup} follow
for all $\lambda$.
\end{proof}

\section{Yangians, weight functions and $R$-matrices}

We recall some facts about the Yangian $Y(\gl_2)$, and explain how to get
weight functions as matrix elements of evaluation Yangian modules.
For more details on Yangians, see for example \cite{CP}.

The Yangian $Y(\gl_2)$ is a Hopf algebra with generators $T_{ij}^{(k)}$,
$i,j = 1,2$; $k= 1,2,3,\dots$. The relations are most conveniently written
in terms of the generating series
$$
T_{ij}(u) = \delta_{ij} + \sum_{k=1}^\infty T_{ij}^{(k)} u^{-k}\,,
$$
and have the form
$$
(u-v) [T_{ij}(u),T_{kl}(v)] = T_{kj}(v)T_{il}(u) - T_{kj}(u)T_{il}(v)\,.
$$
The coproduct for $Y(\gl_2)$ is given by
$$
T_{ij}(u) \mapsto \sum_{k=1}^2 T_{kj}(u) \otimes T_{ik}(u)\,.
$$

The Yangian $Y(\gl_2)$ contains $\mathcal U(\gl_2)$ as a Hopf subalgebra,
the embedding $\mathcal U(\gl_2)\to Y(\gl_2)$ being given by
$E_{ij}\mapsto T_{ji}^{(1)}$. We identify $\mathcal U(\gl_2)$ wih its image
in $Y(\gl_2)$ under this embedding.

There is a family of homomorphisms $\epsilon_z: Y(\gl_2)\to\mathcal U(\gl_2)$,
depending on a complex parameter $z$:
$$
\epsilon_z: T_{ij}(u) \mapsto \delta_{ij} + \frac {E_{ji}}{u-z}\,.
$$
These homomorphism are identical on the subalgebra
$\mathcal U(\gl_2)\subset Y(\gl_2)$.

For any $\gl_2$-module $V$, we denote $V(z)$ the pullback of $V$ through the
homomorphism $\epsilon_z$; the $Y(\gl_2)$-module $V(z)$ is called an evaluation
module.

The Yangian action in a tensor product of evaluation modules is the underlying
structure for the functional realization $\phi[\mathbf l]$ of the tensor
product of the corresponding $\gl_2$-modules.

As before, $\mathbf l = (l_1,l_2,l_3)$ will always denote a triple of
nonnegative integers. For any $k \in \Z_{\ge0}$ set
\begin{gather*}
\mathcal Z_k = \{(k_1,k_2,k_3) \in \Z_{\ge0}^3\,\ |\ \,k_1+k_2+k_3 = k\}\,,
\\[2pt]
\mathcal Z_k[\mathbf l] = \{(k_1,k_2,k_3) \in \mathcal Z_k\,\ |
\ \,k_1\le l_1, k_2 \le l_2, k_3 \le l_3\}.
\end{gather*}

\begin{prop}\cite{KBI,TV3} \label{thm:weightrep}
Let $\mathbf z = (z_1,z_2,z_3)$ be a triple of complex numbers. Consider
the $Y(\gl_2)$-module $V_{(l_1,0)}(z_1) \otimes V_{(l_2,0)}(z_2) \otimes
V_{(l_3,0)}(z_3)$. For any $k \in \Z_{\ge0}, \mathbf q \in \mathcal Z_k[\mathbf
l]$, we have
\begin{align*}
\( \prod_{i=1}^3 \prod_{a=1}^k (t_a-z_i) \) T_{12}(t_1) \dots T_{12}(t_k)
(\1_{l_1} \otimes \1_{l_2} \otimes \1_{l_3}) = \qquad
\\
\sum_{\mathbf p \in \mathcal Z_k[\mathbf l]}
w_{\mathbf p}[\mathbf l; \mathbf z](t_1,\dots,t_k)
E_{12}^{p_1}\1_{l_1}\otimes E_{12}^{p_2}\1_{l_2}\otimes E_{12}^{p_3}\1_{l_3},
\end{align*}
\begin{align*}
\( \prod_{i=1}^3 \prod_{a=1}^k (t_a-z_i) \) T_{21}(t_1) \dots T_{21}(t_k)
(E_{12}^{q_1}\1_{l_1}\otimes E_{12}^{q_2}\1_{l_2}\otimes E_{12}^{q_3}\1_{l_3})
=\qquad
\\
w'_{\mathbf q}[\mathbf l; \mathbf z](t_1,\dots,t_k)
(\1_{l_1} \otimes \1_{l_2} \otimes \1_{l_3})\,,
\end{align*}
where
\begin{align*}
& w_{\mathbf k}[\mathbf l; \mathbf z](t_1,\dots,t_k) =
\\
& \qquad\frac 1{k_1! k_2! k_3!}\;
\Sym\(\prod_{i=1}^3\(\prod_{a=1}^{k_1+\dots +k_{i-1}}
(t_a-z_i+l_i) \prod_{b=k_1+\dots+k_i+1}^{k_1+k_2+k_3} (t_b-z_i)\)
\prod_{a<b} \frac{t_a-t_b-1}{t_a-t_b}\)
\\[12pt]
& w'_{\mathbf k}[\mathbf l; \mathbf z](t_1,\dots,t_k) =
\\
& \qquad\frac 1{k_1! k_2! k_3!}\;
\Sym\(\prod_{i=1}^3\(\prod_{a=1}^{k_1+\dots +k_{i-1}}
(t_a-z_i) \prod_{b=k_1+\dots+k_i+1}^{k_1+k_2+k_3} (t_b-z_i+l_i)\)
\prod_{a<b} \frac{t_a-t_b+1}{t_a-t_b}\).
\end{align*}
\end{prop}

The functions $w_{\mathbf k}[\mathbf l; \mathbf z]$ and $w'_{\mathbf k}[\mathbf
l; \mathbf z]$ are called the weight functions and the dual weight functions,
respectively. They are polynomials in $t_1,\dots,t_k$, and $z_1,z_2,z_3$.

The definition \eqref{eq:phidef} of the the map $\phi[\mathbf l]$ can be
written as
$$
\phi[\mathbf l;\mathbf z] \(E_{21}^{k_1}\1_{l_1}\otimes
E_{21}^{k_2}\1_{l_2}\otimes E_{21}^{k_3}\1_{l_3} \) =
\frac {k_1! l_1!}{(l_1-k_1)!}\,\frac{k_2! l_2!}{(l_2-k_2)!}\,
\frac {k_3! l_3!}{(l_3-k_3)!}\ w'_{\mathbf k}[\mathbf l; \mathbf z]\,.
$$

Let $\Hch_k$ denote the space of symmetric polynomials in variables
$t_1,\dots,t_k$ of degree at most 2 in each $t_1,\dots,t_k$. Set
$$
\Hch = \bigoplus_{k \in \Z_{\ge0}} \Hch_k\,.
$$
It is clear that $\mathcal H_k=\Hch_k\otimes\C[\mathbf z]$ and
$\mathcal H=\Hch\otimes\C[\mathbf z]$.

\begin{thm}
\cite{TV3}
\label{thm:Y2H}
The space $\Hch$ has a structure of a $Y(\gl_2)$-module, depending on
$\mathbf l, \mathbf z$, with the following properties:
\begin{enumerate}
\item
the $\gl_2$-module structure on $\Hch$ coincides with that, induced by
the $\gl_2$-module structure on $\mathcal H$, described in Proposition
\ref{thm:gl2H};
\item
the map $\,\phi[\mathbf l;\mathbf z]:
V_{(l_1,0)}(z_1) \otimes V_{(l_2,0)}(z_2) \otimes V_{(l_3,0)}(z_3)\to\Hch\,$
is a homomorphism of $Y(\gl_2)$-modules.
\end{enumerate}
\end{thm}

We now use the representation theory of the Yangian $Y(\gl_2)$ to prove the
results, announced in Chapter 6. The following theorems go back to \cite{T}.

\begin{thm}\label{thm:yangianreduce}\cite{T}
\begin{enumerate}
\item Let $l \in \Z_{\ge0}$, and let $U$ be an irreducible finite-dimensional
$Y(\gl_2)$-module. Then for generic $z\in\C$, the $Y(\gl_2)$-modules
$U\otimes V_{(l,0)}(z)$ and $V_{(l,0)}(z)\otimes U$ are irreducible and
isomorphic to each other.
\item
Let ${l,l'\in\Z_{\ge0}}$, and let ${z,z'\in\C}$ be such that ${z-l = z'-r}$
for some $r = 0,1, \dots,\allowbreak \min(l,l')-1$. Then the $Y(\gl_2)$-module $V_{(l,0)}(z) \otimes V_{(l',0)}(z')$ contains a unique proper submodule $\widetilde V_{l,l';r}(z)$.
As a $\gl_2$-module, $\widetilde V_{l,l';r}(z)$ is isomorphic to
\begin{equation}
\widetilde V_{l,l';s}(z) \cong \bigoplus_{s=r+1}^{\min(l,l')} V_{(l+l'-s,s)}.
\label{eq:yangsubmodule}
\end{equation}
\end{enumerate}
\end{thm}

\begin{thm} \cite{T}
Let $l,l' \in \Z_{\ge0}$. There exists a unique linear operator $R_{l,l'}(z)
\in \End_\C(V_{(l,0)} \otimes V_{(l',0)})$, rationally depending on a complex
parameter $z$, such that for generic $z,z'$ the operator
$$
\check R_{l,l'}(z-z') = P \circ R_{l,l'}(z-z'):
V_{l}(z) \otimes V_{l'}(z') \to V_{l'}(z') \otimes V_{l}(z)
$$
is an isomorphism of $Y(\gl_2)$-modules, $P$ being the permutation operator:
$P(u\otimes v)=v\otimes u$, and
$$
R_{l,l'}(z) \1_{l} \otimes \1_{l'} = \1_{l} \otimes \1_{l'}\,.
$$
\end{thm}

We call the operator $R_{l,l'}(z)$ the standard rational $R$-matrix.

\begin{proof}[Proof of Proposition \ref{thm:phiproperties}]
Let $\mathbf z_0 = (z_1,z_2,z_3)$ be a generic point from the hyperplane
$z_1-z_2+l_2-r = 0$, with $r = 0,1, \dots, \min(l_1,l_2)-1$. Then
it follows from Theorem \ref{thm:yangianreduce} that the $Y(\gl_2)$-module
$V_{(l_1,0)}(z_1) \otimes V_{(l_2,0)}(z_2) \otimes V_{(l_3,0)}(z_3)$ has
a unique proper submodule $\widetilde V_{l_1,l_2;r}(z_1) \otimes
V_{(l_3,0)}(z_3)$, while Proposition \ref{thm:weightrep} and
formula \eqref{eq:yangsubmodule} imply that $\ker \phi[\mathbf l; \mathbf z_0]$
contains $\widetilde V_{l_1,l_2;r}(z_1) \otimes V_{(l_3,0)}(z_3)$.
Since $\phi[\mathbf l; \mathbf z_0]$ is a non-trivial homomorphism of
$Y(\gl_2)$-modules, the kernel $\ker \phi[\mathbf l; \mathbf z_0]$
must be a proper $Y(\gl_2)$-submodule, which proves the equality
$$
\ker \phi[\mathbf l; \mathbf z_0]=
\widetilde V_{l_1,l_2;r}(z_1) \otimes V_{(l_3,0)}(z_3)\,.
$$

The proof of the second part of Proposition \ref{thm:phiproperties} is similar.
\end{proof}

\begin{proof}[Proof of Proposition \ref{thm:Rmatrix}]
It follows from Proposition \ref{thm:weightrep} that to determine the image of
a vector $v \in V_{l_1} \otimes V_{l_2} \otimes V_{l_3}$, we must act on it by
the element $\( \prod_{i=1}^3 \prod_{a=1}^k (t_a-z_i) \) T_{21}(t_1) \dots
T_{21}(t_k)$, and then use the linear functional $(\1_{l_1}\otimes
\1_{l_2}\otimes \1_{l_3} )^*$. The $Y(\gl_2)$-intertwining property of the
$R$-matrix implies that we have the commutative diagram
\begin{equation}
\begin{CD}
V_{l_2} \otimes V_{l_1} \otimes V_{l_3} @>{\phi[s_1(\mathbf l);s_1(\mathbf z)]}>>
\mathcal H[s_1(\mathbf l); s_1(\mathbf z)] \\
@AA{R_{l_1l_2}(z_1-z_2) \otimes 1}A @| \\
V_{l_1} \otimes V_{l_2} \otimes V_{l_3} @>{\phi[\mathbf l; \mathbf z]}>>
\mathcal H[\mathbf l; \mathbf z]
\\
@VV{1 \otimes R_{l_2l_3}(z_2-z_3)}V @| \\
V_{l_1} \otimes V_{l_3} \otimes V_{l_2} @>{\phi[s_2(\mathbf l);s_2(\mathbf z)]}>>
\mathcal H[s_2(\mathbf l);s_2(\mathbf z)]
\end{CD}\,.
\label{eq:phiHs}
\end{equation}
Now denote
$$
\check{\mathcal R}_{l,l'}(u) =
\frac{\prod_{j=1}^{l'}(u-l+j)}{\prod_{j=1}^l(-u-l'+j)}\ \check R_{l,l'}(u)\,.
$$
The commutative diagram \eqref{eq:phiHs} and the formulas
$$
\frac {\prod_{j=1}^{l_2} (z_1-z_2-l_1+j) }{\prod_{j=1}^{l_1} (z_2-z_1-l_2+j)} =
\frac { \mathcal X[\mathbf l; \mathbf z] }
{ \mathcal X[s_1(\mathbf l); s_1(\mathbf z)] }\,,
$$
$$
\frac {\prod_{j=1}^{l_3} (z_2-z_3-l_2+j) }{\prod_{j=1}^{l_2} (z_3-z_2-l_3+j)} =
\frac { \mathcal X[\mathbf l; \mathbf z] }
{ \mathcal X[s_2(\mathbf l); s_2(\mathbf z)] }\,,
$$
imply that the operators $\check{\mathcal R}_{l,l'}(u)$ satisfy the conditions
of Proposition \ref{thm:Rmatrix}.
\end{proof}

In the remaining part of the Appendix we construct the map
$\mathcal I[\mathbf l]$ explicitly, and describe its singularities.

Let $\mathbf u = (u_1,u_2,u_3)$ be complex variables, and for any $k \in
\Z_{\ge0}$ set $\widetilde{\mathcal H}_k = \mathcal H_k \otimes \C[\mathbf u]$.
Denote also
$$
\Omega_k[\mathbf u; \mathbf z](t_1,\dots,t_k) =
\prod_{i=1}^3\prod_{a=1}^k \frac 1{(t_a-z_i+u_i)(t_a-z_i)}
\prod_{\substack{a,b=1\\a\ne b}}^k \frac{t_a-t_b}{t_a-t_b-1}\;.
$$
For any function $F = F(t_1,\dots,t_k)$, set
$$
\Res_{(\zeta_1,\dots,\zeta_k)} F =
\Res_{t_k=\zeta_k} \( \Res_{t_{k-1}=\zeta_{k-1}}
\( \dots \Res_{t_1=\zeta_1} F(t_1,\dots,t_k) \)\dots\)\,,
$$

Introduce a symmetric $\C(\mathbf u; \mathbf z)$-valued bilinear form
$\<\cdot,\cdot\>$ on each $\widetilde{\mathcal H}_k$, so that for any
$\varphi,\psi\in\widetilde{\mathcal H}_k$ the value $\<\varphi,\psi\>_{\mathbf
u,\mathbf z}$ of the corresponding rational function at generic $\mathbf
u,\mathbf z$ is given by
\begin{equation}
\<\varphi,\psi\>_{\mathbf u,\mathbf z} = \sum_{\mathbf k \in \mathcal Z_k}
\Res_{\mathbf \tau_{\mathbf k}^{-}[\mathbf z]}
\( \varphi[\mathbf u; \mathbf z] \psi[\mathbf u; \mathbf z] \Omega_k[\mathbf u; \mathbf z] \)
\label{eq:residue1},
\end{equation}
where the ``string'' $\mathbf \tau_{\mathbf k}^{-}[\mathbf z]$ is defined by
$$
\mathbf \tau_{\mathbf k}^{-}[\mathbf z] =
(z_1,\dots,z_1-k_1+1,z_2,\dots,z_2-k_2+1,z_3,\dots,z_3-k_3+1)\,.
$$

{\bf Remark.}
The definition above is motivated by the fact (see \cite{TV4}) that the
integral over the imaginary hyperplane $\Re t_1 = \dots = \Re t_k = 0$,
$$
\frac 1{(2\pi i)^k} \idotsint
\varphi[\mathbf u; \mathbf z](t_1,\dots,t_k)
\psi[\mathbf u; \mathbf z](t_1,\dots,t_k)
\Omega_k[\mathbf u; \mathbf z](t_1,\dots,t_k) dt_1 \dots dt_k\,,
$$
converges when $\Re (u_i) \ll \Re (z_i) < 0, \quad i=1,2,3$,
and represents the rational function $\<\varphi,\psi\>$.

Evaluating the integral by residues, and using the analytic continuation,
we get the more convenient algebraic formula \eqref{eq:residue1};
one can show that it is equivalent to another iterated residue formula:
\begin{equation}
\<\varphi,\psi\>_{\mathbf u,\mathbf z} = (-1)^k \sum_{\mathbf k \in \mathcal Z_k}
\Res_{\mathbf \tau^{+}_{\mathbf k}[\mathbf z - \mathbf u]}
\(\varphi[\mathbf u; \mathbf z] \psi[\mathbf u; \mathbf z] \Omega_k[\mathbf u; \mathbf z]\),
\label{eq:residue2}
\end{equation}
where $\mathbf \tau_{\mathbf k}^{+}[\mathbf z - \mathbf u]$ is defined by
$$\mathbf \tau_{\mathbf k}^{+}[\mathbf z - \mathbf u] = (z_1-u_1,\dots,z_1-u_1+k_1-1,z_2-u_2,\dots,z_2-u_2+k_2-1,z_3-u_3,\dots,z_3-u_3+k_3-1).$$

\begin{lem}\label{thm:longpoles}
For any $\varphi,\psi \in \widetilde{\mathcal H}_k$, the only possible
singularities of $\<\varphi,\psi\>_{\mathbf u, \mathbf z}$ are simple poles
at the following hyperplanes:
\begin{equation}
z_i - u_i = z_j - r, \qquad i,j = 1,2,3, \qquad r = 0,1,\dots,k-1.
\label{eq:singularities}
\end{equation}
\end{lem}

\begin{proof}
Since $\varphi,\psi$ are polynomials, possible singularities of
$\<\varphi,\psi\>_{\mathbf u, \mathbf z}$ are determined by $\Omega_k[\mathbf
u; \mathbf z]$. Using explicit formulas, one can check that
$\Res_{\tau^-_{\mathbf k}[\mathbf z]} (\Omega_k[\mathbf u; \mathbf z])$
has only simple poles located at the hyperplanes
$$
z_i-u_i = z_j -r, \qquad z_i = z_j-r, \qquad i,j = 1,2,3,
\qquad r = 0,1,\dots,k-1\,.
$$
Similarly, $\Res_{\tau^+_{\mathbf k}[\mathbf z-\mathbf u]} (\Omega_k[\mathbf u;
\mathbf z])$ has only simple poles located at the hyperplanes
$$
z_i-u_i = z_j -r, \qquad z_i = z_j-u_j-r, \qquad i,j = 1,2,3,
\qquad r = 0,1,\dots,k-1\,.
$$
Since $\<\varphi,\psi\>_{\mathbf u,\mathbf z}$ is given either by formula
\eqref{eq:residue1} or \eqref{eq:residue2}, the function
$\<\varphi,\psi\>_{\mathbf u,\mathbf z}$ can have poles only at the hyperplanes
which belong to both of these lists, which gives precisely list
\eqref{eq:singularities}.
\end{proof}

Let $\widetilde{\mathcal H}_k[\mathbf l]$, $k\in\Z_{\ge0}$, denote the subspace
of functions $\varphi \in \widetilde{\mathcal H}_k$, satisfying the following
conditions for each $i=1,2,3$ such that $l_i < k$:
$$
\varphi(z_i,\dots,z_i-l_i,t_{l_i+2},\dots,t_k) = 0
\quad \text{ at the hyperplane } u_i = l_i\,.
$$
We have two important imbeddings of $\mathcal H_k[\mathbf l]$ into $\widetilde{\mathcal H}_k[\mathbf l]$, given by
$\varphi\mapsto \varphi^-$,
$$\varphi^-[\mathbf u; \mathbf z](t_1,\dots,t_k) = \varphi[\mathbf z](t_1,\dots,t_k),$$
and $\varphi \mapsto \varphi^+$,
$$\varphi^-[\mathbf u; \mathbf z](t_1,\dots,t_k) = \varphi[\mathbf z - \mathbf u + \mathbf l](t_1,\dots,t_k).$$

For $\mathbf p \in \mathcal Z_k$, let $w_{\mathbf p} = w_{\mathbf p}[\mathbf u; \mathbf z]$ be the weight function defined in Proposition \ref{thm:weightrep}, with integers $\mathbf l = (l_1,l_2,l_3)$ replaced by variables $\mathbf u = (u_1,u_2,u_3)$. It is obvious that $w_{\mathbf p} \in \widetilde{\mathcal H}_k$; moreover, if $\mathbf p \in \mathcal Z_k[\mathbf l]$, then $w_{\mathbf p} \in \widetilde{\mathcal H}_k[\mathbf l]$.

\begin{lem}\label{thm:shortpoles}
Let $\mathbf p \in \mathcal Z_k[\mathbf l]$, and $\varphi \in \mathcal H_k[\mathbf l]$.
\begin{enumerate}
\item
The function $\<w_{\mathbf p},\varphi^-\>_{\mathbf u,\mathbf z}$ may have only simple poles, located at the hyperplanes
$$z_i-u_i = z_j-r, \qquad 1 \le i \le j \le 3, \qquad r = 0,1,\dots,l_j-1.$$
\item
The function $\<w_{\mathbf p},\varphi^+\>_{\mathbf u,\mathbf z}$ may have only simple poles, located at the hyperplanes
$$z_i-u_i = z_j-r, \qquad 1 \le i \le j \le 3, \qquad r = 0,1,\dots,l_i-1.$$
\end{enumerate}
\end{lem}

\begin{proof}
Since $\varphi^-[\mathbf u; \mathbf z](\tau^-_{\mathbf q}[\mathbf z]) =
\varphi[\mathbf z](\tau^-_{\mathbf q}[\mathbf z]) = 0$ unless
$\mathbf q \in \mathcal Z_k[\mathbf l]$, we have
\begin{align*}
\<w_{\mathbf p},\varphi^-\>_{\mathbf u, \mathbf z} = {}
& \sum_{\mathbf q \in \mathcal Z_k}\varphi^-(\tau^-_{\mathbf q}[\mathbf z])
\Res_{\tau^-_{\mathbf q}[\mathbf z]}
(w_{\mathbf p}[\mathbf u; \mathbf z] \Omega_k[\mathbf u; \mathbf z]) =
\\
& \sum_{\mathbf q \in \mathcal Z_k[\mathbf l]}
\varphi^-(\tau^-_{\mathbf q}[\mathbf z]) \Res_{\tau^-_{\mathbf q}[\mathbf z]}
(w_{\mathbf p}[\mathbf u; \mathbf z] \Omega_k[\mathbf u; \mathbf z])\,.
\end{align*}

Let $\mathbf q \in \mathcal Z_k[\mathbf l]$. From the explicit formulas for
$w_{\mathbf p}[\mathbf u; \mathbf z]$ and $\Omega_k[\mathbf u; \mathbf z]$ it
follows that the function $\Res_{\tau^-_{\mathbf q}[\mathbf z]} (w_{\mathbf
p}[\mathbf u; \mathbf z] \Omega_k[\mathbf u; \mathbf z])$ has no pole at the
hyperplane $z_i - u_i = z_j - r$ provided $i>j$ or $r \ge q_j$. Using the
description \eqref{eq:singularities} of all possible poles, we see that
$\<w_{\mathbf p},\varphi^-\>_{\mathbf u,\mathbf z}$ may have only poles
corresponding to $i \le j$ and $r < q_j \le l_j$, which proves the first
assertion.

Similarly, using the fact that $\varphi^+[\mathbf u; \mathbf z](\tau^+_{\mathbf
q}[\mathbf z]) \varphi[\mathbf z - \mathbf u](\tau^+_{\mathbf q}[\mathbf z -
\mathbf u]) = 0$ unless $\mathbf q \in \mathcal Z_k[\mathbf l]$, we arrive to
the formula
$$
\<w_{\mathbf p},\varphi^+\>_{\mathbf u, \mathbf z} =
\sum_{\mathbf q \in \mathcal Z_k[l]}
\varphi^+(\tau^+_{\mathbf q}[\mathbf z - \mathbf u])
\Res_{\tau^+_{\mathbf q}[\mathbf z - \mathbf u]}
(w_{\mathbf p}[\mathbf u; \mathbf z] \Omega_k[\mathbf u; \mathbf z])\,.
$$

For any $\mathbf q \in \mathcal Z_k[\mathbf l]$ the function
$\Res_{\tau^+_{\mathbf q}[\mathbf z - \mathbf u]} (w_{\mathbf p}[\mathbf u;
\mathbf z] \Omega_k[\mathbf u; \mathbf z])$ has no pole at the hyperplane $z_i
- u_i = z_j - r$ provided $i>j$ or $r \ge q_i$. Therefore, it follows that
$\<w_{\mathbf p},\varphi^+\>_{\mathbf u; \mathbf z}$ may only have poles
corresponding to $i \le j$ and $r < q_i \le l_i$.
\end{proof}

\begin{lem}\label{thm:analcont}
Let $\varphi,\psi \in \widetilde{\mathcal H}_k[\mathbf l]$.
The function $\<\varphi,\psi\>_{\mathbf u,\mathbf z}$ analytically continues to
$\mathbf u = \mathbf l$, and we have
\begin{align}
\<\varphi,\psi\>_{\mathbf l,\mathbf z} &
= \sum_{\mathbf k \in \mathcal Z_k[\mathbf l]}
\Res_{\mathbf \tau_{\mathbf k}^{-}[\mathbf z]}
\( \varphi[\mathbf l; \mathbf z] \psi[\mathbf l; \mathbf z]
\Omega_k[\mathbf l; \mathbf z] \)
\notag
\\
 & = (-1)^k \sum_{\mathbf k \in \mathcal Z_k[\mathbf l]}
\Res_{\mathbf \tau^{+}_{\mathbf k}[\mathbf z - \mathbf l]}
\(\varphi[\mathbf l; \mathbf z] \psi[\mathbf l; \mathbf z]
\Omega_k[\mathbf l; \mathbf z]\).
\label{eq:analcont}
\end{align}
\end{lem}

In other words, the analytic continuation to $\mathbf u = \mathbf l$ can be
computed by applying the residue formula to the specialized functions
$\psi[\mathbf l;\mathbf z], \varphi[\mathbf l;\mathbf z]$ and $\Omega_k[\mathbf
l;\mathbf z]$, rather than treating $\mathbf u$ as variables, and the summation
is reduced to $\mathbf k \in \mathcal Z_k[\mathbf l] \subset \mathcal Z_k$.
Note that in general this is not true for arbitrary $\varphi,\psi\in
\widetilde{\mathcal H}_k$.

\begin{proof}
One can check that for any $\mathbf q \in \mathcal Z_k$ the rational functions
$\Res_{\tau^-_{\mathbf q}[\mathbf z]} (\psi[\mathbf u; \mathbf z] \Omega_k[\mathbf u; \mathbf z])$
and
$\Res_{\tau^+_{\mathbf q}[\mathbf z - \mathbf u]} (\psi[\mathbf u; \mathbf z] \Omega_k[\mathbf u; \mathbf z])$
have no poles at the hyperplanes $u_i = l_i,\ i = 1,2,3$. Together with Lemma \ref{thm:longpoles}, this shows that the analytic continuation to $\mathbf u = \mathbf l$ is well-defined.

Since for $\mathbf q \notin \mathcal Z_k[\mathbf l]$, the expression $\varphi[\mathbf u; \mathbf z](\tau_{\mathbf k}^-[\mathbf z])$ vanishes at $\mathbf u = \mathbf l$, it follows that the summation in the residue formulas \eqref{eq:residue1} reduces to $\mathbf q \in \mathcal Z_k[\mathbf l]$. Similar conclusion holds for the redisue formula \eqref{eq:residue2}.

The statement now follows from the fact that for $\mathbf q \in \mathcal
Z_k[\mathbf l]$ we have
$$
\Res_{\tau^-_{\mathbf q}[\mathbf z]}(\psi[\mathbf u; \mathbf z]
\Omega_k[\mathbf u; \mathbf z])|_{\mathbf u = \mathbf l} =
\Res_{\tau^-_{\mathbf q}[\mathbf z]} (\psi[\mathbf l; \mathbf z]
\Omega_k[\mathbf l; \mathbf z])\,,
$$
$$\Res_{\tau^+_{\mathbf q}[\mathbf z - \mathbf u]}(\psi[\mathbf u; \mathbf z]
\Omega_k[\mathbf u; \mathbf z])|_{\mathbf u = \mathbf l} =
\Res_{\tau^+_{\mathbf q}[\mathbf z - \mathbf l]} (\psi[\mathbf l; \mathbf z]
\Omega_k[\mathbf l; \mathbf z])\,.
$$
\end{proof}

Introduce a partial order $\le$ on $\mathcal Z_k$ by saying that
$\mathbf x \le \mathbf y$ if and only if $x_1 \le y_1$ and $x_3 \ge y_3$.

\begin{lem}
\label{thm:wwtau}
For any $\mathbf k, \mathbf p,\mathbf q\in\mathcal Z_k$ we have
$$
w_{\mathbf p}[\mathbf l,\mathbf z](\tau_{\mathbf k}^{-}[\mathbf z]) = 0 ,
\qquad \text{ unless }\ \mathbf k \le \mathbf p\,,
$$
$$
w'_{\mathbf q}[\mathbf l,\mathbf z](\tau_{\mathbf k}^{-}[\mathbf z]) = 0 ,
\qquad \text{ unless }\ \mathbf k \ge \mathbf q\,.
$$
Moreover, if \,$\mathbf z$ is generic, then
\,$w'_{\mathbf q}[\mathbf l,\mathbf z](\tau_{\mathbf q}^{-}[\mathbf z])\ne 0$.
\end{lem}
\noindent
The proof is straightforward.

\begin{lem}\cite{TV4}
Let $ k \in \Z_\ge0$, and let $\mathbf p, \mathbf q \in \mathcal Z_k[\mathbf
l]$. Then, we have
\begin{equation}
\<w_{\mathbf p},w'_{\mathbf q}\>_{\mathbf l,\mathbf z} =
\delta_{\mathbf p,\mathbf q}\frac {p_1!l_1!}{(l_1-p_1)!}
\frac {p_2!l_2!}{(l_2-p_2)!}\frac {p_3!l_3!}{(l_3-p_3)!}\;.
\label{eq:dualweightbasis}
\end{equation}
\end{lem}

\begin{proof}
By Lemma \ref{thm:wwtau} the summation in the residue formula
\eqref{eq:analcont} reduces to $\mathbf k$ such that
$\mathbf q \le \mathbf k \le \mathbf p$, that is
$$
\<w_{\mathbf p},w'_{\mathbf q}\>_{\mathbf l,\mathbf z} =
\sum_{\mathbf q \le \mathbf k \le \mathbf p}
w'_{\mathbf q}[\mathbf l; \mathbf z](\tau^-_{\mathbf k}[\mathbf z])
w_{\mathbf p}[\mathbf l; \mathbf z](\tau^-_{\mathbf k}[\mathbf z])
\Res_{\tau^-_{\mathbf k}[\mathbf z]} ( \Omega_k[\mathbf l; \mathbf z])\,.
$$
It follows immediately that $\<w_{\mathbf p},w'_{\mathbf q}\>_{\mathbf
l,\mathbf z}=0$ unless $\mathbf q \le \mathbf p$.

Similarly, one can check that
$$
w_{\mathbf p}(\tau_{\mathbf k}^{+}[\mathbf z-\mathbf l]) = 0 ,
\qquad \text{ unless }\ \mathbf k \ge \mathbf p\,,
$$
$$
w'_{\mathbf q}(\tau_{\mathbf k}^{+}[\mathbf z-\mathbf l]) = 0 ,
\qquad \text{ unless }\ \mathbf k \le \mathbf q\,,
$$
Thus, it follows that
$$
\<w_{\mathbf p},w'_{\mathbf q}\>_{\mathbf l,\mathbf z} =
\sum_{\mathbf p \le \mathbf k \le \mathbf q}
w_{\mathbf p}[\mathbf l; \mathbf z](\tau_{\mathbf k}^{+}[\mathbf z-\mathbf l])
w'_{\mathbf q}[\mathbf l; \mathbf z](\tau_{\mathbf k}^{+}[\mathbf z-\mathbf l])
\( \Res_{\tau_{\mathbf k}^{+}[\mathbf z - \mathbf l]}
\Omega_k[\mathbf l,\mathbf z] \).
$$
and therefore, $\<w_{\mathbf p},w'_{\mathbf q}\>_{\mathbf l,\mathbf z}=0$ unless
$\mathbf p \le \mathbf q$.

Hence, $\<w_{\mathbf p}[\mathbf l,\mathbf z],w'_{\mathbf q}[\mathbf l,\mathbf
z]\>_{\mathbf l,\mathbf z}$ is equal to zero unless $\mathbf p = \mathbf q$.
In the last case there is only one term in the residue summation formula,
corresponding to $\mathbf k = \mathbf p = \mathbf q$, and the statement follows
from the direct computation of the only remaining residue at $\tau_{\mathbf
k}[\mathbf z]$.
\end{proof}

Finally, we have all ingredients necessary to describe explicitly the map
$\mathcal I[\mathbf l]$ and to prove its properties.

\begin{proof}[Proof of Proposition \ref{thm:phibasic}]
Lemma \ref{thm:wwtau} implies that the dual weight functions $\{w'_{\mathbf
k}[\mathbf l; \mathbf z]\}$ are linearly independent as elements of
$\mathcal H_k$, which yields the injectivity of the map $\phi[\mathbf l]$.
Using the explicit combinatorial formulas for $w'_{\mathbf k}[\mathbf l;
\mathbf z]$, one can check that the admissibility conditions \eqref{eq:adm}
are satisfied provided that $\mathbf k\in\mathcal Z_k[\mathbf l]$. Hence,
$\operatorname{Im} \phi[\mathbf l]\otimes_{\C[\mathbf z]} \C(\mathbf z)
\subset \mathcal H_k[\mathbf l]\otimes_{\C[\mathbf z]} \C(\mathbf z)$.
To show that these spaces coincide, it is enough to check that for generic
$\mathbf z$ one has $\dim\mathcal H_k[\mathbf l;\mathbf z]\le
\dim\Im\phi[\mathbf l;\mathbf z]=\#\mathcal Z_k[\mathbf l]$.

Consider the maps $\eps_{\mathbf p}[\mathbf z]:\Hch_k\to\C$,
$$
\eps_{\mathbf p}[\mathbf z]:\psi\mapsto\psi(\tau_{\mathbf p}[\mathbf z])\,.
$$
Lemma \ref{thm:wwtau} implies that for generic $\mathbf z$ these maps are
linearly independent. It is clear that for any $\psi\in\mathcal H[\mathbf
l;\mathbf z]$, $\mathbf p\not\in\mathcal Z_k[\mathbf l]$ one has
$\eps_{\mathbf p}[\mathbf z](\psi)=0$, which gives
$\dim\mathcal H_k[\mathbf l;\mathbf z]\le\#\mathcal Z_k[\mathbf l]$,
since $\dim\Hch_k=\#\mathcal Z_k$.

Now, consider the map
$$
\mathcal I[\mathbf l; \mathbf z]: \mathcal H[\mathbf l; \mathbf z] \to
V_{(l_1,0)} \otimes V_{(l_2,0)} \otimes V_{(l_3,0)}\,,
$$
$$
\mathcal I[\mathbf l; \mathbf z] \varphi = \sum_{k \in \Z_{\ge_0}}
\sum_{\mathbf p \in \mathcal Z_k[\mathbf l]}
\<w_{\mathbf p},\varphi\>_{\mathbf l,\mathbf z}
E_{21}^{p_1} \1_{l_1} \otimes E_{21}^{p_2} \1_{l_2} \otimes
E_{21}^{p_3} \1_{l_3}\,.
$$
Formulas \eqref{eq:phidef} and \eqref{eq:dualweightbasis} imply
\eqref{eq:inversephi}.

According to Lemma \ref{thm:analcont}, we have
$$
\<w_\mathbf p,\varphi\>_{\mathbf l,\mathbf z} =
\(\<w_\mathbf p,\varphi^-\>_{\mathbf u,\mathbf z}\)|_{\mathbf u = \mathbf l}
\(\<w_\mathbf p,\varphi^+\>_{\mathbf u,\mathbf z}\)|_{\mathbf u = \mathbf l}\,,
$$
and Lemma \ref{thm:shortpoles} implies that $\<w_\mathbf p,\varphi\>_{\mathbf
l,\mathbf z}$ may have poles only at the hyperplanes
$$
z_i-u_i = z_j - r,\qquad 1\le i<j\le 3,\quad r = 0,1,\dots,\min(l_i,l_j)-1\,.
$$
This proves the last part of Proposition \ref{thm:phibasic}.
\end{proof}

\section{\strut\kern-.4em}

In this section we prove the existence of the distinguished polynomial basis
$\{\varsigma_i[\mathbf l]\}$ of $\mathcal H_k[\mathbf l]^\sing$, see Lemma
\ref{thm:constantbasis}.

\begin{proof}[Proof of Lemma \ref{thm:constantbasis}.]
The subspace $\mathcal H_k[\mathbf l]^\sing$ consists of polynomials of degree
at most 1 in each $t_j$, and therefore any function $\varsigma(t_1,\dots,t_k)
\in \mathcal H_k[\mathbf l]^\sing$ has the form
$$
\varsigma(t_1,\dots,t_k)=
\sum_{r=0}^k X_r(z_1,z_2,z_3)\,\sigma_r(t_1,\dots,t_k)
$$
where $X_0,\ldots X_k$ are some coefficients from $\C[\mathbf z]$, and
$$
\sigma_r(t_1,\dots,t_k) =
\sum_{1\le i_1<\dots<i_r\le k}t_{i_1}\dots t_{i_r}, \qquad r = 0,1,2,\dots;
$$
in particular, $\sigma_0(t_1,\dots,t_k)\equiv 1$,
and $\sigma_r(t_1,\dots,t_k) \equiv 0$ if $r>k$.

Using the combinatorial identity
$$
\sigma_r(z,\dots,z-l,t_{l+2},\dots,t_k) =
\sum_{s=0}^r \sigma_{r-s}(z,\dots,z-l) \sigma_s(t_{l+2},\dots,t_k)\,,
$$
we see that admissibility conditions \eqref{eq:adm} for the function
$\varsigma(t_1,\ldots,t_k)$ are equivalent to the following system of equations
for the coefficients $X_0,\ldots X_k$:
\begin{equation}
\sum_{r=s}^k \,X_r\,\sigma_{r-s}(z_i,\dots,z_i-l_i) = 0,
\qquad s = 0,\dots,k-l_i-1, \quad i =1,2,3.
\label{eq:Xsystem}
\end{equation}

Let ${d = \dim \(V_{(l_1,0)} \otimes V_{(l_2,0)} \otimes
V_{(l_3,0)}\)^\sing[(l_1+l_2+l_3-k,k)]}$. We will show that system
\eqref{eq:Xsystem} determines $X_0,\dots,X_{k-d}$ as linear combinations of
free variables $X_{k-d+1},\dots,X_k$ with coefficients, polynomial in
$z_1,z_2,z_3$. Then we construct $\varsigma_s[\mathbf l],\quad s = 1,\dots,d$,
by taking
$$
X_{k-s+1}=1,\qquad X_{k-r+1}=0\quad\text{ for }\ r\ne s, \quad 1\le r \le d\,,
$$
and computing $X_0,\dots,X_{k-d}$ from equations \eqref{eq:Xsystem}.
The polynomials $\varsigma_s[\mathbf l]$ obeys invariance \eqref{eq:invaradm}
by construction.

Consider the $(k-d+1) \times (k+1)$ matrix of system \eqref{eq:Xsystem}.
It suffices to prove that its $(k-d+1) \times (k-d+1)$ minor, corresponding
to the variables $X_0,\dots,X_{k-d}$, is nonzero, and divides all other
$(k-d+1)\times (k-d+1)$ minors. Before giving a general proof of this
statement, we consider an example.

{\bf Example:} ${l_1 = 1}$, ${l_2 = 2}$, ${k = 3}$, ${l_3 > 3}$. Then
${d=\dim\(V_{(1,0)}\otimes V_{(2,0)}\otimes V_{(l_3,0)}\)^\sing[(l_3,3)]}
\allowbreak=1$, and the matrix of system of \eqref{eq:Xsystem} is equal to
$$
\begin{pmatrix}
1 & z_1 + (z_1-1) & z_1(z_1-1) & 0 \\
0 & 1 & z_1 + (z_1-1) & z_1(z_1-1) \\
1\, & z_2+(z_2-1)+(z_2-2)\, & z_2(z_2-1) + z_2(z_2-2) + (z_2-1)(z_2-2)\,&
z_2(z_2-1)(z_2-2)
\end{pmatrix}.
$$
The $3\times3$ minors of this matrix are given by
$$
M_{012} = 3(z_1-z_2)(z_1-z_2+1)\,,
$$
$$
M_{013} = (z_1-z_2)(z_1-z_2+1)(2z_1+z_2-2)\,,
$$
$$
M_{023} = (z_1-z_2)(z_1-z_2+1)(z_1^2 + 2 z_1 z_2 - 3z_1 - z_2 + 2)\,,
$$
$$
M_{123}=(z_1-z_2)(z_1-z_2+1)(3 z_1^2z_2 -3 z_1^2 -3 z_1 z_2 + 5z_1 + z_2 -2)\,.
$$
The leftmost minor $M_{012}$ is nonzero and divides all other minors.
The function $\varsigma(t_1,t_2,t_3)$, spanning the space
$\mathcal H[(1,2,l_3),(z_1,z_2,z_3)]^\sing[(l_3,3)]$, is given by
$$
\varsigma(t_1,t_2,t_3) = t_1t_2t_3 -
\frac1{M_{012}}\,\bigl(M_{013}(t_1t_2 + t_2t_3 + t_1t_3) -
M_{023}(t_1+t_2 +t_3) + M_{123}\bigr)\,.
$$

\goodbreak

We now return to the case of arbitrary $\mathbf l$. Without loss of generality,
we may assume that $l_1 \le l_2 \le l_3$. To shorten notation, we denote
$$
\sigma_r^{(i)} = \sigma_r(z_i,z_i-1,\dots,z_i-l_i),
\qquad r = 0,\dots,k, \quad i = 1,2,3.
$$

{\bf Case 1: $k\le l_1$.} Then, there are no admissibility
conditions, $d = k+1$, and all the variables $X_0,\ldots X_k$ are free.

{\bf Case 2: $l_1 < k \le l_2$.} Then $d = l_1+1$, and the matrix of system
\eqref{eq:Xsystem} is equal to
$$
\begin{pmatrix}
\sigma_{0}^{(1)} && \sigma_{1}^{(1)} && \sigma_{2}^{(1)} && \dots && \sigma_{k-l_1-1}^{(1)} && \sigma_{k-l_1}^{(1)} && \sigma_{k-l_1+1}^{(1)} && \dots && \sigma_k^{(1)}\\
0 && \sigma_{0}^{(1)} && \sigma_{1}^{(1)} && \dots &&\sigma_{k-l_1-2}^{(1)} && \sigma_{k-l_1-1}^{(1)} && \sigma_{k-l_1}^{(1)} && \dots && \sigma_{k-1}^{(1)}\\
0 && 0 && \sigma_{0}^{(1)} && \dots &&\sigma_{k-l_1-3}^{(1)} && \sigma_{k-l_1-2}^{(1)} && \sigma_{k-l_1-1}^{(1)} && \dots && \sigma_{k-2}^{(1)}\\
\vdots && \vdots && \vdots && \ddots && \dots && \dots && \dots && \ddots && \dots\\
0 && 0 && 0 && \dots && \sigma_{0}^{(1)} && \sigma_{1}^{(1)} && \sigma_{2}^{(1)} && \dots && \sigma_{l_1+1}^{(1)}
\end{pmatrix}
$$

Note that some of the terms in the upper right corner can also vanish,
because $\sigma_{l_1+2}^{(1)} = \dots = \sigma_{k}^{(1)} = 0$.

Since $\sigma_0^{(1)} \equiv 1$, the leftmost $(k-d+1) \times (k-d+1)$ minor is
identically equal to 1. Therefore $X_{k-d+1},\dots,X_k$ are free variables,
and $X_0,\dots,X_{k-d}$ are their linear combinations with coefficients,
polynomial in $z_1,z_2,z_3$.

{\bf Case 3: $l_2 < k \le l_3$.} If $k\ge l_1+l_2$, then $d = 0$, and there is
nothing to check. Assume that $k<l_1+l_2$, which gives $d = l_1+l_2-k$.
The matrix of system \eqref{eq:Xsystem} is equal to
\begin{equation}
\begin{pmatrix}
\sigma_{0}^{(1)} && \sigma_{1}^{(1)} && \dots && \dots && \dots &&\sigma_{k-l_1-1}^{(1)} && \dots && \sigma_k^{(1)}\\
0 && \sigma_{0}^{(1)} && \dots && \dots && \dots &&\sigma_{k-l_1-2}^{(1)} && \dots && \sigma_{k-1}^{(1)}\\
\vdots && \vdots && \vdots && \ddots && \dots && \dots && \ddots && \dots\\
0 && 0 && \dots && \dots && \dots && \sigma_{0}^{(1)} && \dots && \sigma_{l_1+1}^{(1)} \\
\sigma_{0}^{(2)} && \sigma_{1}^{(2)} && \dots && \sigma_{k-l_2-1}^{(2)} && \dots && \dots && \dots && \sigma_k^{(2)}\\
0 && \sigma_{0}^{(2)} && \dots &&\sigma_{k-l_2-2}^{(2)} && \dots && \dots && \dots && \sigma_{k-1}^{(2)}\\
\vdots && \vdots && \ddots && \dots && \dots && \dots && \ddots && \dots\\
0 && 0 && \dots && \sigma_{0}^{(2)} && \dots && \dots && \dots && \sigma_{l_2+1}^{(2)}
\end{pmatrix}
\label{eq:matrix}
\end{equation}

We make two claims now. First, every $(k-d+1)\times(k-d+1)$ minor of this
matrix is a polynomial in $z_1,z_2$, divisible by
$$
D(z_1,z_2) =
\prod_{j_1=1}^{k-l_1} \prod_{j_2=1}^{k-l_2} (z_1-z_2+j_1-j_2)\,.
$$
Second, the leftmost $(k-d+1)\times(k-d+1)$ minor of matrix \eqref{eq:matrix}
is a nonzero polynomial in $z_1,z_2$ of homogeneous degree $(k-l_1)(k-l_2)$.

These two claims together imply the desired statement. Indeed, since
the homogeneous degree of the polynomial $D(z_1,z_2)$ is equal to
$(k-l_1)(k-l_2)$, the leftmost minor is proportional to $D(z_1,z_2)$ with
a nonzero coefficient. Therefore, every $(k-d+1)\times(k-d+1)$ minor is
divisible by the leftmost $(k-d+1)\times(k-d+1)$ minor, which means that
$X_0,\dots,X_{k-d}$ can be expressed as linear combinations of the free
variables $X_{k-d+1},\dots,X_k$ with coefficients, polynomially depending on
$z_1,z_2,z_3$.

We first prove the second claim. Note that $\sigma_r^{(i)}$ is a polynomial
of $z_i$ with the highest degree term $\binom {l_i+1}r z_i^r$. Therefore,
the homogeneous degree of the leftmost $(k-d+1)\times (k-d+1)$ minor does not
exceed $(k-l_1)(k-l_2)$, while the sum of terms of degree $(k-l_1)(k-l_2)$
is given by the determinant
$$
\det\begin{pmatrix}
\vphantom{\Big|}
1 & \binom{l_1+1}1 z_1 & \dots & \dots & \dots &\binom{l_1+1}{k-l_1-1} z_1^{k-l_1-1} & \dots & \binom {l_1+1}{k-d-1} z_1^{k-d-1}
\vphantom{\Big|}\\
0 & 1 & \dots & \dots & \dots &\binom{l_1+1}{k-l_1-2} z_1^{k-l_1-2} & \dots & \binom{l_1+1}{k-d-2} z_1^{k-d-2}
\vphantom{\Big|}\\
\vdots & \vdots & \vdots & \ddots & \dots & \dots & \ddots & \dots
\vphantom{\Big|}\\
0 & 0 & \dots & \dots & \dots & 1 & \dots & \binom{l_1+1}{k-l_2} z_1^{k-l_2}
\vphantom{\Big|}\\
1 & \binom{l_2+1}1 z_2 & \dots & \binom{l_2+1}{k-l_2-1}z_2^{k-l_2-1} & \dots & \dots & \dots & \binom{l_2+1}{k-d-1} z_2^{k-d-1}
\vphantom{\Big|}\\
0 & 1 & \dots &\binom{l_2+1}{k-l_2-2} z_2^{k-l_2-2} & \dots & \dots & \dots & \binom{l_2+1}{k-d-2} z_2^{k-d-2}
\vphantom{\Big|}\\
\vdots & \vdots & \ddots & \dots & \dots & \dots & \ddots & \dots
\vphantom{\Big|}\\
0 & 0 & \dots & 1 & \dots & \dots & \dots & \binom{l_2+1}{k-l_1} z_2^{k-l_1}
\end{pmatrix}
$$
The monomial $z_1^{(k-l_1)(k-l_2)}$ enters this determinant with a nonzero
coefficient
$$
\det\begin{pmatrix}
\binom{l_1+1}{k-l_2} & \binom{l_1+1}{k-l_2+1} &\dots & \binom {l_1+1}{2k-l_1-l_2-1}
\vphantom{\Big|}\\
\binom{l_1+1}{k-l_2-1} & \binom{l_1+1}{k-l_2} & \dots & \binom{l_1+1}{k-1}
\vphantom{\Big|}\\
\vdots & \vdots & \ddots & \vdots
\vphantom{\Big|}\\
\binom{l_1+1}{l_1-l_2+1} & \binom{l_1+1}{l_1-l_2+2} & \dots & \binom{l_1+1}{k-l_2}
\end{pmatrix}\,=\;
\prod_{i=1}^{k-l_1}\,\prod_{j=1}^{k-l_2}\,\prod_{s=1}^{l_1+l_2-k+1}\,
\frac{i+j+s-1}{i+j+s-2}\;.
$$
The last equality follows from the general formula, established in \cite{DET}:
$$
\det\begin{pmatrix}
\binom{a+b}{a} & \binom{a+b}{a+1} & \dots & \binom {a+b}{a+c}
\vphantom{\Big|}\\
\binom{a+b}{a-1} & \binom{a+b}{a} & \dots & \binom {a+b}{a+c-1}
\vphantom{\Big|}\\
\vdots & \vdots & \ddots & \vdots
\vphantom{\Big|}\\
\binom{a+b}{a-c} & \binom{a+b}{a-c+1} & \dots & \binom{a+b}{a}
\end{pmatrix}\,=\;
\prod_{i=1}^{c+1}\,\prod_{j=1}^a\,\prod_{s=1}^b\;\frac{i+j+s-1}{i+j+s-2}\;,
$$
by taking $a=k-l_2$, $b=l_1+l_2-k+1$, $c= k-l_1-1$.

Finally, we prove the first claim. Recall that we assume
$l_1\le l_2<k<l_1+l_2$.

Take an integer $j$ such that $-l_2\le j\le l_1$ and
let $z_1-z_2-j=0$. Then the sets $\{z_1,z_1-1,\dots,z_1-l_1\}$ and
$\{z_2,z_2-1,\dots,z_2-l_2\}$ have $p(j)+1$ common points
$\{\zeta,\zeta-1,\dots,\zeta-p(j)\}$, where $\zeta = z_1$ for $j<0$,
$\zeta=z_2$ for $j\ge 0$, and
$$
p(j) = \left\{
\begin{aligned}
& l_2+j, & \qquad -l_2 \le j < l_1-l_2\,, \\
& l_1, & \qquad l_1-l_2 \le j < 0\,, \\
& l_1-j, & \qquad 0 \le j \le l_1\,,\\
\end{aligned}
\right.
$$

Consider the system of linear equations on variables $X_0,\dots,X_k$,
corresponding to the condition
\begin{equation}
\varsigma(\zeta,\zeta-1,\dots,\zeta-p(j),t_{p(j)+2},\dots,t_k) \equiv 0\,.
\label{eq:zeta}
\end{equation}
Let $L\subset \C^k$ denote the linear span of the rows of the corresponding
matrix. Since the rows are linearly independent, we have $\dim L = k-p(j)$.

Obviously, if the function $\varsigma(t_1,\dots,t_k)$ obeys condition
\eqref{eq:zeta}, then
\begin{gather*}
\varsigma(z_1,z_1-1,\dots,z_1-l_1,t_{l_1+1},\dots,t_k) \equiv 0\,,
\\
\varsigma(z_2,z_2-1,\dots,z_2-l_2,t_{l_2+1},\dots,t_k) \equiv 0\,.
\end{gather*}
This means that the linear span of the rows of matrix \eqref{eq:matrix} is
contained in $L$. In particular, the rank of matrix \eqref{eq:matrix} does not
exceed $\dim L$. Therefore, if $l_1-k<j<k-l_2$, then the rank of matrix
\eqref{eq:matrix} is less than its number of rows $2k-l_1-l_2$ by at least
$k+p(j)-l_1-l_2$, which is positive.

The above consideration shows that every $(k-d+1)\times(k-d+1)$ minor of
matrix \eqref{eq:matrix} is divisible by the product
$$
\prod_{j=l_1-k+1}^{k-l_2-1}(z_1-z_2-j)^{k+p(j)-l_1-l_2}\,.
$$
It is easy to show that this product coincides with $D(z_1,z_2)$.

{\bf Case 4: $l_1 + l_2 < k \le l_3$.}
In this case $\dim \mathcal H_k(\mathbf l;\mathbf z)^\sing = 0$, and there is nothing to check.

{\bf Case 5: $k>l_3$.}
There are three types of admissibility conditions in this case, related to
$(z_1,l_1)$, $(z_2,l_2)$ and $(z_3,l_3)$. The consideration, similar to Case 3,
shows that the leftmost $(k-d+1)\times(k-d+1)$ minor of the corresponding
matrix is nonzero and divides every $(k-d+1)\times(k-d+1)$ minor.
This implies the required statement.
\end{proof}

\end{document}